\documentclass{gtart}
\def\ifplaintex{\expandafter\ifx\csname documentclass\endcsname\relax}

\def\gtp{{\mathsurround=0pt\it $\cal G\mskip-2mu$eometry \&\ 
$\cal T\!\!$opology $\cal P\!$ublications}}  % GT publications

\def\recd{{\small Received:\qua\receiveddate\ifx\reviseddate\relax
\else\qquad Revised:\qua\reviseddate\fi\par}} 

\def\lognumber#1{\def\thelognumber{#1}}
\def\volumenumber#1{\def\thevolumenumber{#1}}
\def\volumeyear#1{\def\thevolumeyear{#1}}
\def\papernumber#1{\def\thepapernumber{#1}}
\def\pagenumbers#1#2{\def\startpage{#1}\def\finishpage{#2}}
\def\published#1{\def\publishdate{#1}}

\def\received#1{\def\receiveddate{#1}}
\def\revised#1{\def\reviseddate{#1}}
\def\accepted#1{\def\accepteddate{#1}}

\long\def\asciiabstract#1{\long\def\theasciiabstract{#1}}

\let\\\par\let\thelognumber\relax\let\thevolumenumber\relax
\let\thepapernumber\relax\let\thevolumeyear\relax\let\startpage\relax
\let\finishpage\relax\let\publishdate\relax\let\receiveddate\relax
\let\reviseddate\relax\let\accepteddate\relax\let\theasciititle\relax
\let\theasciiauthors\relax
\let\theasciiabstract\relax

\let\theasciiemail\relax

\ifplaintex
\font\logobig=cmssbx10 scaled 3836
\font\logomed=cmssbx10 scaled 2557
\else
\font\logobig=cmssbx10 scaled 4200
\font\logomed=cmssbx10 scaled 2800
\fi

\long\def\makeagttitle{   %%% start of definition of \makeagttitle
\count0=\startpage
\agt\hfill      %   Journal title (top left) 
\hbox to 45truept{\vbox to 0pt{\vglue -13truept{\logomed A\kern -.37em{\logobig 
T}\kern -.38em G}\vss}\hss}
\break
{\small Volume \thevolumenumber\ (\thevolumeyear)
\startpage--\finishpage\nl
Published: \publishdate}

\vglue .25truein

{\parskip=0pt\leftskip 0pt plus
1fil\def\\{\par\smallskip}{\Large\bf\thetitle}\par\medskip} \vglue
0.05truein

{\parskip=0pt\leftskip 0pt plus 1fil\def\\{\par}{\sc\theauthors}
\par\medskip}%
 
\vglue 0.03truein 

{\small\leftskip 25truept\rightskip 25truept{\bf Abstract}\stdspace\theabstract

{\bf AMS Classification}\stdspace\theprimaryclass
\ifx\thesecondaryclass\relax\else; \thesecondaryclass\fi\par
{\bf Keywords}\stdspace \thekeywords\par}\vglue 7truept

}   %%%% end of definition of \makeagttitle

\ifplaintex
\font\phead=cmsl9 scaled 950
\font\pnum=cmbx10 scaled 913
\font\pfoot=cmsl9 scaled 950
\headline{\vbox to 0pt{\vskip -4.5mm\line{\small\phead\ifnum
\count0=\startpage ISSN 1472-2739 (on-line) 1472-2747 (printed)
\hfill {\pnum\folio}\else\ifodd\count0\def\\{ }% 
\ifx\theshorttitle\relax\thetitle\else\theshorttitle\fi\hfill{\pnum\folio}
\else\def\\{ and }{\pnum\folio}\hfill\ifx\theshortauthors\relax\theauthors
\else\theshortauthors\fi\fi\fi}\vss}}
\footline{\vbox to 0pt{\vglue 0mm\line{\small\pfoot\ifnum\count0=\startpage
\copyright\ \gtp\hfill\else
\agt, Volume \thevolumenumber\ (\thevolumeyear)\hfill\fi}\vss}}
\else
\font=cmsl9 scaled 1050
\font\lnum=cmbx10 
\font=cmsl9 scaled 1050
\makeatletter
\def\@oddhead{{\small\lhead\ifnum\count0=\startpage ISSN 1472-2739 
(on-line) 1472-2747 (printed)\hfill {\lnum\number\count0}\else\ifodd\count0
\def\\{ }\ifx\theshorttitle\relax \thetitle \else\theshorttitle\fi\hfill
{\lnum\number\count0}\else\def\\{ and }{\lnum\number\count0}
\hfill\ifx\theshortauthors\relax 
\theauthors\else\theshortauthors\fi\fi\fi}}\def\@evenhead{\@oddhead}
\def\@oddfoot{\small\lfoot\ifnum\count0=\startpage\copyright\ \gtp\hfill\else
\agt, Volume \thevolumenumber\ (\thevolumeyear)\hfill\fi}
\def\@evenfoot{\@oddfoot}
\makeatother
\fi
\let\maketitlepage\makeagttitle
\let\makeshorttitle\maketitlepage
\let\maketitle\maketitlepage

\newwrite\gtoutfile
\long\gdef\makeheadfile{  %%% start of definition of \makeheadfile
{\def\\{, }\def\s{ }
\immediate\openout\gtoutfile head.xxx
\immediate\write\gtoutfile{To: math@arxiv.org}
\immediate\write\gtoutfile{Subject: put OR rep NNNNN:ppppp}
\immediate\write\gtoutfile{--text follows this line--}
\immediate\write\gtoutfile{Proxy-for: \ifx\theasciiauthors\relax
\theauthors\else\theasciiauthors\fi\s<\ifx\theasciiemail\relax\theemail\else\theasciiemail\fi>}
\immediate\write\gtoutfile{\noexpand\\}
\immediate\write\gtoutfile{Authors: \ifx\theasciiauthors\relax
\theauthors\else\theasciiauthors\fi}
{\def\\{ }\immediate\write\gtoutfile{Title: \ifx\theasciititle\relax
\thetitle\else\theasciititle\fi}}
\immediate\write\gtoutfile{Subj-class: GT or SG, GR etc}
\immediate\write\gtoutfile{MSC-class: \theprimaryclass\ifx\thesecondaryclass\relax\else, \thesecondaryclass\fi}
\immediate\write\gtoutfile{Journal-ref: Algebr. Geom. Topol. \thevolumenumber\s
(\thevolumeyear) \startpage-\finishpage}
\immediate\write\gtoutfile{Comments: Published by Algebraic and
Geometric Topology at}
\immediate\write\gtoutfile{\s\s\s  http://www.maths.warwick.ac.uk/agt/AGTVol\thevolumenumber/agt-\thevolumenumber-\thepapernumber.abs.html}
\immediate\write\gtoutfile{\noexpand\\}
\immediate\write\gtoutfile{}
\ifx\theasciiabstract\relax
\immediate\write\gtoutfile{\theabstract}\else
\immediate\write\gtoutfile{\theasciiabstract}\fi
\immediate\write\gtoutfile{}
\immediate\write\gtoutfile{\noexpand\\}
\immediate\write\gtoutfile{}
\immediate\closeout\gtoutfile}}  %%% end of definition of \makeheadfile

\def\maketitlepage{\makeagttitle\makeheadfile}
\let\makeshorttitle\maketitlepage
\let\maketitle\maketitlepage

\def\ifplaintex{\expandafter\ifx\csname documentclass\endcsname\relax}

\def\gtp{{\mathsurround=0pt\it $\cal G\mskip-2mu$eometry \&\ 
$\cal T\!\!$opology $\cal P\!$ublications}}  % GT publications

\def\recd{{\small Received:\qua\receiveddate\ifx\reviseddate\relax
\else\qquad Revised:\qua\reviseddate\fi\par}} 

\def\lognumber#1{\def\thelognumber{#1}}
\def\volumenumber#1{\def\thevolumenumber{#1}}
\def\volumeyear#1{\def\thevolumeyear{#1}}
\def\papernumber#1{\def\thepapernumber{#1}}
\def\pagenumbers#1#2{\def\startpage{#1}\def\finishpage{#2}}
\def\published#1{\def\publishdate{#1}}

\def\received#1{\def\receiveddate{#1}}
\def\revised#1{\def\reviseddate{#1}}
\def\accepted#1{\def\accepteddate{#1}}

\long\def\asciiabstract#1{\long\def\theasciiabstract{#1}}

\let\\\par\let\thelognumber\relax\let\thevolumenumber\relax
\let\thepapernumber\relax\let\thevolumeyear\relax\let\startpage\relax
\let\finishpage\relax\let\publishdate\relax\let\receiveddate\relax
\let\reviseddate\relax\let\accepteddate\relax\let\theasciititle\relax
\let\theasciiauthors\relax
\let\theasciiabstract\relax

\let\theasciiemail\relax

\ifplaintex
\font\logobig=cmssbx10 scaled 3836
\font\logomed=cmssbx10 scaled 2557
\else
\font\logobig=cmssbx10 scaled 4200
\font\logomed=cmssbx10 scaled 2800
\fi

\long\def\makeagttitle{   %%% start of definition of \makeagttitle
\count0=\startpage
\agt\hfill      %   Journal title (top left) 
\hbox to 45truept{\vbox to 0pt{\vglue -13truept{\logomed A\kern -.37em{\logobig 
T}\kern -.38em G}\vss}\hss}
\break
{\small Volume \thevolumenumber\ (\thevolumeyear)
\startpage--\finishpage\nl
Published: \publishdate}

\vglue .25truein

{\parskip=0pt\leftskip 0pt plus
1fil\def\\{\par\smallskip}{\Large\bf\thetitle}\par\medskip} \vglue
0.05truein

{\parskip=0pt\leftskip 0pt plus 1fil\def\\{\par}{\sc\theauthors}
\par\medskip}%
 
\vglue 0.03truein 

{\small\leftskip 25truept\rightskip 25truept{\bf Abstract}\stdspace\theabstract

{\bf AMS Classification}\stdspace\theprimaryclass
\ifx\thesecondaryclass\relax\else; \thesecondaryclass\fi\par
{\bf Keywords}\stdspace \thekeywords\par}\vglue 7truept

}   %%%% end of definition of \makeagttitle

\ifplaintex
\font\phead=cmsl9 scaled 950
\font\pnum=cmbx10 scaled 913
\font\pfoot=cmsl9 scaled 950
\headline{\vbox to 0pt{\vskip -4.5mm\line{\small\phead\ifnum
\count0=\startpage ISSN 1472-2739 (on-line) 1472-2747 (printed)
\hfill {\pnum\folio}\else\ifodd\count0\def\\{ }% 
\ifx\theshorttitle\relax\thetitle\else\theshorttitle\fi\hfill{\pnum\folio}
\else\def\\{ and }{\pnum\folio}\hfill\ifx\theshortauthors\relax\theauthors
\else\theshortauthors\fi\fi\fi}\vss}}
\footline{\vbox to 0pt{\vglue 0mm\line{\small\pfoot\ifnum\count0=\startpage
\copyright\ \gtp\hfill\else
\agt, Volume \thevolumenumber\ (\thevolumeyear)\hfill\fi}\vss}}
\else
\font=cmsl9 scaled 1050
\font\lnum=cmbx10 
\font=cmsl9 scaled 1050
\makeatletter
\def\@oddhead{{\small\lhead\ifnum\count0=\startpage ISSN 1472-2739 
(on-line) 1472-2747 (printed)\hfill {\lnum\number\count0}\else\ifodd\count0
\def\\{ }\ifx\theshorttitle\relax \thetitle \else\theshorttitle\fi\hfill
{\lnum\number\count0}\else\def\\{ and }{\lnum\number\count0}
\hfill\ifx\theshortauthors\relax 
\theauthors\else\theshortauthors\fi\fi\fi}}\def\@evenhead{\@oddhead}
\def\@oddfoot{\small\lfoot\ifnum\count0=\startpage\copyright\ \gtp\hfill\else
\agt, Volume \thevolumenumber\ (\thevolumeyear)\hfill\fi}
\def\@evenfoot{\@oddfoot}
\makeatother
\fi
\let\maketitlepage\makeagttitle
\let\makeshorttitle\maketitlepage
\let\maketitle\maketitlepage

\newwrite\gtoutfile
\long\gdef\makeheadfile{  %%% start of definition of \makeheadfile
{\def\\{, }\def\s{ }
\immediate\openout\gtoutfile head.xxx
\immediate\write\gtoutfile{To: math@arxiv.org}
\immediate\write\gtoutfile{Subject: put OR rep NNNNN:ppppp}
\immediate\write\gtoutfile{--text follows this line--}
\immediate\write\gtoutfile{Proxy-for: \ifx\theasciiauthors\relax
\theauthors\else\theasciiauthors\fi\s<\ifx\theasciiemail\relax\theemail\else\theasciiemail\fi>}
\immediate\write\gtoutfile{\noexpand\\}
\immediate\write\gtoutfile{Authors: \ifx\theasciiauthors\relax
\theauthors\else\theasciiauthors\fi}
{\def\\{ }\immediate\write\gtoutfile{Title: \ifx\theasciititle\relax
\thetitle\else\theasciititle\fi}}
\immediate\write\gtoutfile{Subj-class: GT or SG, GR etc}
\immediate\write\gtoutfile{MSC-class: \theprimaryclass\ifx\thesecondaryclass\relax\else, \thesecondaryclass\fi}
\immediate\write\gtoutfile{Journal-ref: Algebr. Geom. Topol. \thevolumenumber\s
(\thevolumeyear) \startpage-\finishpage}
\immediate\write\gtoutfile{Comments: Published by Algebraic and
Geometric Topology at}
\immediate\write\gtoutfile{\s\s\s  http://www.maths.warwick.ac.uk/agt/AGTVol\thevolumenumber/agt-\thevolumenumber-\thepapernumber.abs.html}
\immediate\write\gtoutfile{\noexpand\\}
\immediate\write\gtoutfile{}
\ifx\theasciiabstract\relax
\immediate\write\gtoutfile{\theabstract}\else
\immediate\write\gtoutfile{\theasciiabstract}\fi
\immediate\write\gtoutfile{}
\immediate\write\gtoutfile{\noexpand\\}
\immediate\write\gtoutfile{}
\immediate\closeout\gtoutfile}}  %%% end of definition of \makeheadfile

\def\maketitlepage{\makeagttitle\makeheadfile}
\let\makeshorttitle\maketitlepage
\let\maketitle\maketitlepage

\def\ifplaintex{\expandafter\ifx\csname documentclass\endcsname\relax}

\def\gtp{{\mathsurround=0pt\it $\cal G\mskip-2mu$eometry \&\ 
$\cal T\!\!$opology $\cal P\!$ublications}}  % GT publications

\def\recd{{\small Received:\qua\receiveddate\ifx\reviseddate\relax
\else\qquad Revised:\qua\reviseddate\fi\par}} 

\def\lognumber#1{\def\thelognumber{#1}}
\def\volumenumber#1{\def\thevolumenumber{#1}}
\def\volumeyear#1{\def\thevolumeyear{#1}}
\def\papernumber#1{\def\thepapernumber{#1}}
\def\pagenumbers#1#2{\def\startpage{#1}\def\finishpage{#2}}
\def\published#1{\def\publishdate{#1}}

\def\received#1{\def\receiveddate{#1}}
\def\revised#1{\def\reviseddate{#1}}
\def\accepted#1{\def\accepteddate{#1}}

\long\def\asciiabstract#1{\long\def\theasciiabstract{#1}}

\let\\\par\let\thelognumber\relax\let\thevolumenumber\relax
\let\thepapernumber\relax\let\thevolumeyear\relax\let\startpage\relax
\let\finishpage\relax\let\publishdate\relax\let\receiveddate\relax
\let\reviseddate\relax\let\accepteddate\relax\let\theasciititle\relax
\let\theasciiauthors\relax
\let\theasciiabstract\relax

\let\theasciiemail\relax

\ifplaintex
\font\logobig=cmssbx10 scaled 3836
\font\logomed=cmssbx10 scaled 2557
\else
\font\logobig=cmssbx10 scaled 4200
\font\logomed=cmssbx10 scaled 2800
\fi

\long\def\makeagttitle{   %%% start of definition of \makeagttitle
\count0=\startpage
\agt\hfill      %   Journal title (top left) 
\hbox to 45truept{\vbox to 0pt{\vglue -13truept{\logomed A\kern -.37em{\logobig 
T}\kern -.38em G}\vss}\hss}
\break
{\small Volume \thevolumenumber\ (\thevolumeyear)
\startpage--\finishpage\nl
Published: \publishdate}

\vglue .25truein

{\parskip=0pt\leftskip 0pt plus
1fil\def\\{\par\smallskip}{\Large\bf\thetitle}\par\medskip} \vglue
0.05truein

{\parskip=0pt\leftskip 0pt plus 1fil\def\\{\par}{\sc\theauthors}
\par\medskip}%
 
\vglue 0.03truein 

{\small\leftskip 25truept\rightskip 25truept{\bf Abstract}\stdspace\theabstract

{\bf AMS Classification}\stdspace\theprimaryclass
\ifx\thesecondaryclass\relax\else; \thesecondaryclass\fi\par
{\bf Keywords}\stdspace \thekeywords\par}\vglue 7truept

}   %%%% end of definition of \makeagttitle

\ifplaintex
\font\phead=cmsl9 scaled 950
\font\pnum=cmbx10 scaled 913
\font\pfoot=cmsl9 scaled 950
\headline{\vbox to 0pt{\vskip -4.5mm\line{\small\phead\ifnum
\count0=\startpage ISSN 1472-2739 (on-line) 1472-2747 (printed)
\hfill {\pnum\folio}\else\ifodd\count0\def\\{ }% 
\ifx\theshorttitle\relax\thetitle\else\theshorttitle\fi\hfill{\pnum\folio}
\else\def\\{ and }{\pnum\folio}\hfill\ifx\theshortauthors\relax\theauthors
\else\theshortauthors\fi\fi\fi}\vss}}
\footline{\vbox to 0pt{\vglue 0mm\line{\small\pfoot\ifnum\count0=\startpage
\copyright\ \gtp\hfill\else
\agt, Volume \thevolumenumber\ (\thevolumeyear)\hfill\fi}\vss}}
\else
\font=cmsl9 scaled 1050
\font\lnum=cmbx10 
\font=cmsl9 scaled 1050
\makeatletter
\def\@oddhead{{\small\lhead\ifnum\count0=\startpage ISSN 1472-2739 
(on-line) 1472-2747 (printed)\hfill {\lnum\number\count0}\else\ifodd\count0
\def\\{ }\ifx\theshorttitle\relax \thetitle \else\theshorttitle\fi\hfill
{\lnum\number\count0}\else\def\\{ and }{\lnum\number\count0}
\hfill\ifx\theshortauthors\relax 
\theauthors\else\theshortauthors\fi\fi\fi}}\def\@evenhead{\@oddhead}
\def\@oddfoot{\small\lfoot\ifnum\count0=\startpage\copyright\ \gtp\hfill\else
\agt, Volume \thevolumenumber\ (\thevolumeyear)\hfill\fi}
\def\@evenfoot{\@oddfoot}
\makeatother
\fi
\let\maketitlepage\makeagttitle
\let\makeshorttitle\maketitlepage
\let\maketitle\maketitlepage

\newwrite\gtoutfile
\long\gdef\makeheadfile{  %%% start of definition of \makeheadfile
{\def\\{, }\def\s{ }
\immediate\openout\gtoutfile head.xxx
\immediate\write\gtoutfile{To: math@arxiv.org}
\immediate\write\gtoutfile{Subject: put OR rep NNNNN:ppppp}
\immediate\write\gtoutfile{--text follows this line--}
\immediate\write\gtoutfile{Proxy-for: \ifx\theasciiauthors\relax
\theauthors\else\theasciiauthors\fi\s<\ifx\theasciiemail\relax\theemail\else\theasciiemail\fi>}
\immediate\write\gtoutfile{\noexpand\\}
\immediate\write\gtoutfile{Authors: \ifx\theasciiauthors\relax
\theauthors\else\theasciiauthors\fi}
{\def\\{ }\immediate\write\gtoutfile{Title: \ifx\theasciititle\relax
\thetitle\else\theasciititle\fi}}
\immediate\write\gtoutfile{Subj-class: GT or SG, GR etc}
\immediate\write\gtoutfile{MSC-class: \theprimaryclass\ifx\thesecondaryclass\relax\else, \thesecondaryclass\fi}
\immediate\write\gtoutfile{Journal-ref: Algebr. Geom. Topol. \thevolumenumber\s
(\thevolumeyear) \startpage-\finishpage}
\immediate\write\gtoutfile{Comments: Published by Algebraic and
Geometric Topology at}
\immediate\write\gtoutfile{\s\s\s  http://www.maths.warwick.ac.uk/agt/AGTVol\thevolumenumber/agt-\thevolumenumber-\thepapernumber.abs.html}
\immediate\write\gtoutfile{\noexpand\\}
\immediate\write\gtoutfile{}
\ifx\theasciiabstract\relax
\immediate\write\gtoutfile{\theabstract}\else
\immediate\write\gtoutfile{\theasciiabstract}\fi
\immediate\write\gtoutfile{}
\immediate\write\gtoutfile{\noexpand\\}
\immediate\write\gtoutfile{}
\immediate\closeout\gtoutfile}}  %%% end of definition of \makeheadfile

\def\maketitlepage{\makeagttitle\makeheadfile}
\let\makeshorttitle\maketitlepage
\let\maketitle\maketitlepage

\def\ifplaintex{\expandafter\ifx\csname documentclass\endcsname\relax}

\def\gtp{{\mathsurround=0pt\it $\cal G\mskip-2mu$eometry \&\ 
$\cal T\!\!$opology $\cal P\!$ublications}}  % GT publications

\def\recd{{\small Received:\qua\receiveddate\ifx\reviseddate\relax
\else\qquad Revised:\qua\reviseddate\fi\par}} 

\def\lognumber#1{\def\thelognumber{#1}}
\def\volumenumber#1{\def\thevolumenumber{#1}}
\def\volumeyear#1{\def\thevolumeyear{#1}}
\def\papernumber#1{\def\thepapernumber{#1}}
\def\pagenumbers#1#2{\def\startpage{#1}\def\finishpage{#2}}
\def\published#1{\def\publishdate{#1}}

\def\received#1{\def\receiveddate{#1}}
\def\revised#1{\def\reviseddate{#1}}
\def\accepted#1{\def\accepteddate{#1}}

\long\def\asciiabstract#1{\long\def\theasciiabstract{#1}}

\let\\\par\let\thelognumber\relax\let\thevolumenumber\relax
\let\thepapernumber\relax\let\thevolumeyear\relax\let\startpage\relax
\let\finishpage\relax\let\publishdate\relax\let\receiveddate\relax
\let\reviseddate\relax\let\accepteddate\relax\let\theasciititle\relax
\let\theasciiauthors\relax
\let\theasciiabstract\relax

\let\theasciiemail\relax

\ifplaintex
\font\logobig=cmssbx10 scaled 3836
\font\logomed=cmssbx10 scaled 2557
\else
\font\logobig=cmssbx10 scaled 4200
\font\logomed=cmssbx10 scaled 2800
\fi

\long\def\makeagttitle{   %%% start of definition of \makeagttitle
\count0=\startpage
\agt\hfill      %   Journal title (top left) 
\hbox to 45truept{\vbox to 0pt{\vglue -13truept{\logomed A\kern -.37em{\logobig 
T}\kern -.38em G}\vss}\hss}
\break
{\small Volume \thevolumenumber\ (\thevolumeyear)
\startpage--\finishpage\nl
Published: \publishdate}

\vglue .25truein

{\parskip=0pt\leftskip 0pt plus
1fil\def\\{\par\smallskip}{\Large\bf\thetitle}\par\medskip} \vglue
0.05truein

{\parskip=0pt\leftskip 0pt plus 1fil\def\\{\par}{\sc\theauthors}
\par\medskip}%
 
\vglue 0.03truein 

{\small\leftskip 25truept\rightskip 25truept{\bf Abstract}\stdspace\theabstract

{\bf AMS Classification}\stdspace\theprimaryclass
\ifx\thesecondaryclass\relax\else; \thesecondaryclass\fi\par
{\bf Keywords}\stdspace \thekeywords\par}\vglue 7truept

}   %%%% end of definition of \makeagttitle

\ifplaintex
\font\phead=cmsl9 scaled 950
\font\pnum=cmbx10 scaled 913
\font\pfoot=cmsl9 scaled 950
\headline{\vbox to 0pt{\vskip -4.5mm\line{\small\phead\ifnum
\count0=\startpage ISSN 1472-2739 (on-line) 1472-2747 (printed)
\hfill {\pnum\folio}\else\ifodd\count0\def\\{ }% 
\ifx\theshorttitle\relax\thetitle\else\theshorttitle\fi\hfill{\pnum\folio}
\else\def\\{ and }{\pnum\folio}\hfill\ifx\theshortauthors\relax\theauthors
\else\theshortauthors\fi\fi\fi}\vss}}
\footline{\vbox to 0pt{\vglue 0mm\line{\small\pfoot\ifnum\count0=\startpage
\copyright\ \gtp\hfill\else
\agt, Volume \thevolumenumber\ (\thevolumeyear)\hfill\fi}\vss}}
\else
\font=cmsl9 scaled 1050
\font\lnum=cmbx10 
\font=cmsl9 scaled 1050
\makeatletter
\def\@oddhead{{\small\lhead\ifnum\count0=\startpage ISSN 1472-2739 
(on-line) 1472-2747 (printed)\hfill {\lnum\number\count0}\else\ifodd\count0
\def\\{ }\ifx\theshorttitle\relax \thetitle \else\theshorttitle\fi\hfill
{\lnum\number\count0}\else\def\\{ and }{\lnum\number\count0}
\hfill\ifx\theshortauthors\relax 
\theauthors\else\theshortauthors\fi\fi\fi}}\def\@evenhead{\@oddhead}
\def\@oddfoot{\small\lfoot\ifnum\count0=\startpage\copyright\ \gtp\hfill\else
\agt, Volume \thevolumenumber\ (\thevolumeyear)\hfill\fi}
\def\@evenfoot{\@oddfoot}
\makeatother
\fi
\let\maketitlepage\makeagttitle
\let\makeshorttitle\maketitlepage
\let\maketitle\maketitlepage

\newwrite\gtoutfile
\long\gdef\makeheadfile{  %%% start of definition of \makeheadfile
{\def\\{, }\def\s{ }
\immediate\openout\gtoutfile head.xxx
\immediate\write\gtoutfile{To: math@arxiv.org}
\immediate\write\gtoutfile{Subject: put OR rep NNNNN:ppppp}
\immediate\write\gtoutfile{--text follows this line--}
\immediate\write\gtoutfile{Proxy-for: \ifx\theasciiauthors\relax
\theauthors\else\theasciiauthors\fi\s<\ifx\theasciiemail\relax\theemail\else\theasciiemail\fi>}
\immediate\write\gtoutfile{\noexpand\\}
\immediate\write\gtoutfile{Authors: \ifx\theasciiauthors\relax
\theauthors\else\theasciiauthors\fi}
{\def\\{ }\immediate\write\gtoutfile{Title: \ifx\theasciititle\relax
\thetitle\else\theasciititle\fi}}
\immediate\write\gtoutfile{Subj-class: GT or SG, GR etc}
\immediate\write\gtoutfile{MSC-class: \theprimaryclass\ifx\thesecondaryclass\relax\else, \thesecondaryclass\fi}
\immediate\write\gtoutfile{Journal-ref: Algebr. Geom. Topol. \thevolumenumber\s
(\thevolumeyear) \startpage-\finishpage}
\immediate\write\gtoutfile{Comments: Published by Algebraic and
Geometric Topology at}
\immediate\write\gtoutfile{\s\s\s  http://www.maths.warwick.ac.uk/agt/AGTVol\thevolumenumber/agt-\thevolumenumber-\thepapernumber.abs.html}
\immediate\write\gtoutfile{\noexpand\\}
\immediate\write\gtoutfile{}
\ifx\theasciiabstract\relax
\immediate\write\gtoutfile{\theabstract}\else
\immediate\write\gtoutfile{\theasciiabstract}\fi
\immediate\write\gtoutfile{}
\immediate\write\gtoutfile{\noexpand\\}
\immediate\write\gtoutfile{}
\immediate\closeout\gtoutfile}}  %%% end of definition of \makeheadfile

\def\maketitlepage{\makeagttitle\makeheadfile}
\let\makeshorttitle\maketitlepage
\let\maketitle\maketitlepage

\lognumber{45}
\volumenumber{2}
\volumeyear{2002}
\papernumber{45}
\pagenumbers{1155}{1178}
\received{12 June 2002}
\revised{8 November 2002}
\accepted{17 December 2002}
\published{27 December 2002}

 \usepackage{amsmath}
 \usepackage{amssymb}
 \usepackage{graphicx}

\let\endpf\endproof 

\theoremstyle{definition}
 \newtheorem{defn}{Definition}[section]
 \newtheorem{definition}[defn]{Definition}

 \newtheorem{remark}[defn]{Remark}

\theoremstyle{plain}
 \newtheorem{lemma}[defn]{Lemma}
 
 \newtheorem{theorem}[defn]{Theorem}
 \newtheorem{prop}[defn]{Proposition}
 \newtheorem{corollary}[defn]{Corollary}

 \newcommand{\calf}{{\cal F}}

 \newcommand{\calm}{{\cal M}}

 \newcommand{\sltz}{\textrm{SL}(2, \mathbb{Z})}

 \newcommand{\g}{{\gamma}}
 \newcommand{\G}{{\Gamma}}
 \newcommand{\e}{{\epsilon}}

 \newcommand{\s}{{\sigma}}
\def\e{\epsilon}
\def\p{^\prime}

\def\i{^{-1}}
\def\l{\lambda}

\def\g{\gamma}
\def\s{\sigma}

\def\G{\Gamma}

\def\<{\langle}
\def\>{\rangle}

\begin{document}

 \title{Groups generated by positive multi-twists\\and the fake lantern
problem}
\shorttitle{Groups generated by positive multi-twists}

 \authors{Hessam Hamidi-Tehrani}
 \address{ B.C.C. of the City University of New York\\Bronx, NY
 10453, USA}
\email{hessam@math.columbia.edu}

\begin{abstract}
Let $\G$ be a group generated by two positive multi-twists. We
give some sufficient conditions for $\G$ to be free or have no
``unexpectedly reducible'' elements. For a group $\G$ generated by
two Dehn twists, we classify the elements in $\G$ which are
multi-twists. As a consequence we are able to list all the
lantern-like relations in the mapping class groups. We classify
groups generated by powers of two Dehn twists which are free,  or
have no ``unexpectedly reducible" elements. In the end we pose
similar problems for groups generated by powers of $n \ge 3$
twists and give a partial result.
 \end{abstract}

\asciiabstract{Let Gamma be a group generated by two positive
multi-twists. We give some sufficient conditions for Gamma to be free
or have no `unexpectedly reducible' elements. For a group Gamma
generated by two Dehn twists, we classify the elements in Gamma which
are multi-twists. As a consequence we are able to list all the
lantern-like relations in the mapping class groups. We classify groups
generated by powers of two Dehn twists which are free, or have no
``unexpectedly reducible" elements. In the end we pose similar
problems for groups generated by powers of n > 2 twists and give a
partial result.}

\primaryclass{57M07}                

\secondaryclass{20F38, 57N05}

\keywords {Mapping class group, Dehn twist, multi-twist,
pseudo-Anosov, lantern relation }

\maketitle

\section{Introduction}\label{intro}

In a meeting of the American Mathematical Society in Ann Arbor, MI
in March 2002, John McCarthy posed the following question: Suppose
a collection of simple closed curves satisfy the lantern relation
(see Figure~\ref{s040}) algebraically. Is it true that they must
form a lantern, as in the same figure, given the same
commutativity conditions? In this article we consider groups that
are generated by two multi-twists and give conditions that
guarantee the group is free or does not contain an accidental
multi-twist. This will, in particular, answer McCarthy's question
to the affirmative (see Theorem~\ref{lanternlike}).

 To make
this precise, let $S$ be an oriented surface, possibly with
punctures. For the isotopy class of\footnote{We will usually drop
the phrase ``isotopy class of'' in the rest of this paper for
 brevity, as all curves are considered up to isotopy.} a simple  closed curve $c$
on $S$ let $T_c$ denote the right-handed Dehn twist about $c$.
Let $(c_1,c_2)$ denote the minimum geometric intersection number
 of isotopy classes of 1-sub-manifolds $c_1,c_2.$ By $\calm(S)$ we denote the
 mapping class group of $S$, i.e, the group of homeomorphisms of $S$ which permute the
  punctures, up to isotopies fixing the punctures.

The free group on $n$ generators will be denoted by  $\mathbb
F_n$.

Let $A=\{a_1,\cdots,a_k\}$ be a collection of non-parallel,
non-trivial, pairwise disjoint simple closed curves. For any integers
$m_1,...,m_k$, we call $T_A=T_{a_k}^{m_k} \cdots T_{a_1}^{m_1}$
a multi-twist. If, furthermore, all $m_i>0$, we call $T_A$ a
positive multi-twist. We will study the group generated by two
positive multi-twists in detail. We will give explicit conditions
which imply $\<T_A, T_B\> \cong \mathbb{F}_2$ (see
Theorem~\ref{2multifree}). For a group $\< T_a, T_b \>$ generated
by two Dehn twists, we give a complete description of elements: We
determine which elements are multi-twists, and which elements are
pseudo-Anosov restricted to the subsurface which is a regular
neighborhood of $a \cup b$ (see Theorems~\ref{2twistsfree},
\ref{pa(a,b)=2}, and \ref{(a,b)=2}). A mapping class $f$ is called
pseudo-Anosov if $f^n(c) \ne c$ (up to isotopy) for all
non-trivial simple closed curves $c$ and $n>0$. Let $A=\{a_1,a_2,
..., a_n\}$ be a set of non-parallel, non-trivial simple closed
curves on $S$. The surface filled by $A$, denoted by $S_A$, is a
regular neighborhood $N$ of $a_1 \cup \cdots  \cup a_n$ together
with the components of $S \setminus N$ which are discs with 0 or
 1 puncture, assuming that $a_i$'s are drawn as geodesics
 of some constant-curvature metric. $S_A$ is well-defined independent
 of chosen metric \cite{hs}.  We say that $A$ fills up $S$ if $S_A=S$.

\begin{definition}\label{relativepa}
{\rm{ A word $w=T_{c_1}^{n_1} \cdots T_{c_k}^{n_k}$ is called a
cyclically-reduced word if $c_1 \ne c_k$. For such a word $w$,
define $supp(w)=S_{\{c_1,\cdots, c_k\}}$. Then we say $w$ is {\it{
relatively pseudo-Anosov}} if the restriction of the map $w$ is
pseudo-Anosov in $\calm(U)$, for all components $U$ of $supp(w)$
which are not annuli. If $g=h w h\i$ (as words) with $w$ cyclically-reduced,
define $supp(g)=h(supp(w))$. Then define $g$ to be relatively
pseudo-Anosov in the same way as above. }}
\end{definition}

In the above, the equation $g=hwh\i$ is an equation of words, not elements, otherwise one can easily give examples where the definition breaks down. To show that the above definition is well-defined, note that
if $w=T_{c_1}^{n_1} \cdots T_{c_k}^{n_k}$ is such that $c_1=c_k$
but $n_1 \ne -n_k$, then one can write $w=T_{c_1}^{n_1} w'
T_{c_1}^{-n_1}=T_{c_1}^{-n_k} w'' T_{c_1}^{n_k}$, where $w',w''$
are both cyclically reduced. Notice that
$T_{c_1}^n(S_{\{c_1,\cdots,c_{k-1}\}})=S_{\{c_1,\cdots,c_{k-1}\}}$
for all $n$ and so $supp(w)=\{c_1,\cdots, c_{k-1}\}$.

Also note that a power of a Dehn twist $T_{c_i}^{n_i}$ is
a relatively pseudo-Anosov word since its support is an annulus.
Similarly a multi-twist is relatively pseudo-Anosov as well.
 A group $\G$ with given set of multi-twist generators is relatively pseudo-Anosov if every reduced word in
generators of $\G$ is relatively pseudo-Anosov.

 Intuitively, a group $\G$ generated by multi-twists is relatively
 pseudo-Anosov if no word in generators of $\G$
has ``unexpected reducibility''.

It should be noticed that in the case of two curves
 $a,b$ filling up a closed surface this was done by Thurston as a method
 to construct pseudo-Anosov elements;
i.e., he showed that $\<T_a,T_b\>$ is free and consists of
pseudo-Anosov elements besides powers of conjugates of the
generators \cite{FLP}. Our methods are completely different and
elementary, and are only based on how the geometric intersection
behaves under Dehn twists.

One surprising result that we find is a lantern-like relation:
$$(T_bT_a)^2=T_{\partial_1}T_{\partial_2}T_{\g}^{-4}T_{\g'}^{-4},$$
where these curves are defined in Figures~\ref{abins120} and
\ref{gamma} (see Proposition~\ref{tbta2}). This relation is
lantern-like in the sense that the left hand side is a word in two
Dehn twists about intersecting curves and the right hand side is a
multi-twist. We then prove that this relation and the lantern
relation are the only lantern-like relations
(Theorem~\ref{lanternlike}).

In the case when $n \ge 3$, we give some sufficient conditions for
$\G=\<T_{a_1},\cdots, T_{a_n} \>$ to be isomorphic to $F_n$. To motivate
our condition, look at the case
 $\G=\<T_{a_1}, T_{a_2}, T_{a_3} \>$, and assume
$a_3=T_{a_1}(a_2)$. Now $T_{a_3}=T_{a_1}T_{a_2}T_{a_1}\i$, so $\G
\ncong \mathbb F_3$. But notice that  $(a_1,a_3)=(a_1,a_2)$ and
$(a_2,a_3)=(a_1,a_2)^2$, by Lemma~\ref{flplemma}. This shows that
the set $I=\{(a_i,a_j) \ | \ i\ne j \}$ is ``spread around''. It
turns out that this is in a sense an obstruction for $\G \cong \mathbb
F_n$:

\medskip{\bf Theorem}\qua {\sl Suppose $\G= \<
T_{a_1},..., T_{a_h}\>$, and let $m=\min I$ and $M=\max I$,
where $I=\{(a_i,a_j) \ | \ i\ne j \}$. Then $\G \cong \mathbb F_h$
if $M \le m^2/6$.}

\medskip

We will prove a more general version of this (see
Theorem~\ref{3.2}).

It should be noticed that similar arguments have been used to
prove that certain groups generated by three $2 \times 2$ matrices
are free \cite{BM, Sch}.

In Section~\ref{basics} we go over basic facts about Dehn twists
and geometric intersection pairing and different kinds of
ping-pong arguments we are going to use. In Section~\ref{2dt} we
prove our general theorems about groups generated by two positive
multi-twists. In Section~\ref{S20} we look at the specific case of
a lantern formation. In Section~\ref{S12} we look at a formation
which produces a lantern-like relation. In Section~\ref{applan} we
prove that the only possible lantern-like relations are the ones
given in Theorem~\ref{lantern}. In Section~\ref{nge3} we prove a
theorem on groups generated by $n$ Dehn twists. In Section 8 we
pose some questions that are of similar flavor.

\begin{remark} {\rm

 After the completion of this work, the author learned that Dan
Margalit has obtained some results on the subject of lantern
relation  using the action of the mapping class group on homology
\cite{M}. Also notice that Theorem~\ref{lantern} here answers the
first question in \cite[Section 7]{M}. }
\end{remark}

\section{Basics}\label{basics}

 For two isotopy classes of closed 1-sub-manifolds  $a,b$ of $S$ let $(a,b)$ denote their geometric
intersection number. For a set of closed 1-sub-manifolds
$A=\{a_1,..., a_n\}$ and a simple closed curve $x$ put
$$||x||_A= \sum_{i=1}^n (x,a_i).$$
For a non-trivial simple closed curve let $T_a$ be the
(right-handed) Dehn twist in curve $a$. The following lemma is
proved in \cite{FLP}.

\begin{lemma}\label{flplemma} For simple closed curves $a,x,b$, and $n \ge 0$,
$$ |(T^{\pm n}_a (x),b)-n(x,a)(a,b)| \le (x,b).$$
\end{lemma}

Let $a=\{a_1,..., a_k\}$ be a collection of distinct, mutually disjoint non-trivial
isotopy classes of simple closed curves. For integers $n_i>0$, the mapping class $T_a=T^{n_1}_{a_1}\cdots T^{n_k}_{a_k}$ is called a positive multi-twist. We also have the following lemma:

\begin{lemma}\label{flplemma2} For a positive multi-twist $T_a=T^{n_1}_{a_1}\cdots
T^{n_k}_{a_k}$, 1-sub-manifolds $x,b$ and $n \in \mathbb{Z}$,
$$ |(T^{n}_a (x),b)-|n|\sum_{i=1}^k n_i(x,a_i)(a_i,b)| \le (x,b).$$
\end{lemma}

For a proof see \cite[Lemma 4.2]{I}. The statement of that lemma
has the expression $|n|-2$ instead of $|n|$ above. Using the
assumption that all $n_i$ are positive, the same proof goes
through to prove the improved statement given here. Alternatively,
a proof can be found in \cite[Expos\'e 4]{FLP}.

The classic ping-pong argument was used first by Klein \cite{K}.
We give two versions here which will be applied in
Section~\ref{2dt}. The group $\G$ can be a general group. The notation
 $\G=\<f_1,..., f_n\>$ means that the group $\G$ is generated by elements
$f_1,...,f_n$.

\begin{lemma}[Ping-pong]\label{pp}Let $\G=\<f_1,..., f_n\>$, $n\ge 2$. Suppose
 $\G$ acts on a set $X$. Assume that
 there are $n$  non-empty mutually disjoint subsets $X_1,\cdots,X_n$ of $X$
such that $f_i^{\pm k}(\cup_{j \ne i} X_j) \subset X_i$, for all
$1 \le i \le n$ and $k>0$. Then $\G \cong \mathbb F_n$.

\end{lemma}

\proof First notice that a non-empty reduced  word of form $w=f_1^*f_i^*
\cdots f_j^* f_1^*$ (*'s are non-zero integers) is not the
identity because $w(X_2) \cap X_2 \subset X_1 \cap X_2=
\emptyset$. But any reduced word in $f_1^{\pm 1},\cdots,f_n^{\pm
1}$ is conjugate to a $w$ of the above form. \endpf

\begin{lemma}[Tower ping-pong]\label{tpp}Let $\G$ be a
group generated by $f_1,\cdots, f_n$. Suppose
 $\G$ acts on a set $X$, and there is a function
$||.||:X \to \mathbb R_{\ge 0}$, with the following properties:
 There are $n$  non-empty mutually disjoint subsets $X_1,\cdots,X_n$ of $X$
such that $f_i^{\pm k}(X \setminus X_i) \subset X_i$ and for any
$x \in X \setminus X_i$, we have $||f_i^{\pm k}(x)||>||x||$  for
all $k>0$. Then $\G \cong \mathbb F_n$. Moreover, the action of
every $g \in \G$ which is not conjugate to some power of some
$f_i$   on $X$ has no periodic points.
\end{lemma}

\proof Any non-empty reduced  word in $f_1^*,...,f_n^*$ (*'s denote
non-zero integers) is conjugate to a reduced word $w=f_1^*\cdots
f_1^*$. To show that $w \ne id$ notice that if $x_1 \in X
\setminus X_1$, then $w(x_1) \in X_1$, therefore $w(x_1) \ne x_1$.
To prove the last assertion, notice that it's enough to show the
claim with ``periodic points'' replaced by ``fixed points''. Any
element of $\G$ which is not conjugate to a power of some $f_i$ is
conjugate to some reduced word of the form $w=f_j^* \cdots f_i^*$
with $i \ne j$. Now suppose $w(x)=x$. First assume $x \in X
\setminus X_i$. Then by assumption $||w(x)||>||x||$ which is
impossible. If on the other hand, $x \in X_i$ and $w(x)=x$, then
$w\i(x)=f_i^* \cdots f_j^*(x)=x$. But again by assumption
$||w\i(x)||>||x||$, which is a contradiction. \endpf

\renewcommand\labelenumi{\rm(\roman{enumi})}

\section{Groups generated by two positive multi-twists}\label{2dt}

Let $a=\{a_1,\cdots,a_k\}$ and $b=\{b_1,\cdots, b_l\}$ be two collections
 of isotopy classes of non-trivial, mutually disjoint simple closed curves
  on $S$, respectively, such that $(a,b)>0$. Let $m_1,\cdots, m_k,$ and $n_1,\cdots,n_l$ be
   positive integers. In this section we will set

$\bullet$\qua $T_a=T_{a_1}^{m_1}\cdots T_{a_k}^{m_k}$ and $T_b=T_{b_1}^{n_1}\cdots T_{b_l}^{n_l}$.

$\bullet$\qua $A=\{a, b \}$.

$\bullet$\qua $X=\{x \ | \ x {\textrm{ is the isotopy class of a simple closed
curve and }} ||x||_A >0 \ \}.$

$\bullet$\qua For $\l \in (0,\infty)$ set
\begin{gather*}
N_{a,\l}=\{x \in X \ | \  \ (x,a)<\l (x,b)\},\\
N_{b, \l\i}=\{x \in X\ | \  \ \l (x,b)<(x,a)\}.
\end{gather*}
Notice that $a \in N_{a,\l}$ and $b \in N_{b, \l\i}$,
 and $N_{a,\l} \cap N_{b,\l\i}=\emptyset$.
Moreover, $\<T_a,T_b\>$ acts on $X$, and when $\l$ is irrational $X=N_{a,\l} \cup N_{b,\l\i}$.
\begin{lemma}\label{sets}
With the above notation:

\begin{enumerate}
\item $T_a^{\pm n}(N_{b,\l\i}) \subset N_{a,\l}$ if
 $n m_i(a_i,b)\ge 2\l\i$ for all $1 \le i \le k$.
\item If $n m_i(a_i,b) \ge 2\l\i$ for all $1 \le i \le k$,
and $x \in  N_{b,\l\i}$, then $|| T_a^{\pm n}(x)||_A > ||x||_A.$
\item $T_b^{\pm n}(N_{a,\l})
\subset N_{b,\l\i}$ if $n n_j (a,b_j)\ge 2\l$ for all $1 \le j \le l$.
\item If $n n_j (a, b_j) \ge 2\l$ for all $1 \le j \le l$, and $x
\in N_{a,\l}$, then $|| T_b^{\pm n}(x)||_A > ||x||_A.$
\end{enumerate}
\end{lemma}

\proof Suppose $x \in N_{b, \l\i}$, and $n>0$ such that $n m_i
(a_i,b)\ge 2\l\i$. Then by Lemma~\ref{flplemma2},
\begin{eqnarray*} (T^{\pm n}_a(x), b) &\ge& n \sum_i m_i (x,a_i)(a_i,b)-(x,b) \\
            &>& n\sum_i m_i(x, a_i)(a_i,b)-\l\i\sum_i (x,a_i) \\
        &=& \sum_i (n m_i(a_i,b)-\l\i) (x,a_i) \\
            &\ge& \l\i \sum_i (x,a_i) \\
            &=& \l\i (x,a) \\
        &=& \l\i (T^{\pm n}_a(x), T^{\pm n}_a(a)) \\
        &=& \l\i (T^{\pm n}_a(x), a).\\
\end{eqnarray*}
This proves (i).
 By symmetry we immediately get (iii). Now notice that for $x \in N_{b, \l\i}$,
 \begin{eqnarray*} ||T^{\pm n }_a(x)||_A &=&(T^{\pm n }_a(x),a)+(T^{\pm n }_a(x),b) \\
            &\ge& (x,a)+ n \sum_i m_i(x,a_i)(a_i,b)-(x,b) \\
            &>& \sum_i (1+n m_i(a_i,b)-\l\i)(x,a_i) \\
            &=& \sum_i\l(1+n m_i(a_i,b)-\l\i)(1+\l)\i(\l\i(x,a_i)+(x,a_i)). \\
            \end{eqnarray*}
 But $\l(1+n m_i(a_i,b)-\l\i)(1+\l)\i \ge 1$ if and only if  $n m_i(a_i,b) \ge 2\l\i$, which by assumption implies
\begin{eqnarray*}
 ||T^{\pm n }_a(x)||_A &>&\sum_i(\l\i(x,a_i)+(x,a_i)) \\
                &=&\l\i(x,a)+(x,a)\\
                &>& (x,b)+(x,a) \\
            &=&  ||x||_A.
\end{eqnarray*}

This proves (ii), and by symmetry (iv).\endpf

\begin{theorem}\label{2multifree} For two positive multi-twists
 $T_a=T^{m_1}_{a_1}\cdots T^{m_k}_{a_k}$ and $T_b=T^{n_1}_{b_1}\cdots T^{n_l}_{b_l}$ on the surface $S$, the group $\< T_a, T_b\> \cong \mathbb F_2$ if both of the following conditions are satisfied:

 \begin{enumerate}
 \item $m_i(a_i,b)\ge 2$ for all $1 \le i \le k$.
 \item $n_j(a,b_j)\ge 2$ for all $1 \le j \le l$.
\end{enumerate}
\end{theorem}

\proof The group $\<T_a, T_b \>$ acts on $X=\{x \ | \ ||x||_A
>0\}$, where $A=\{a, b\}$. Now use the sets $X_1=N_{a,1}$ and $X_2=N_{b,
1}$ in Lemma~\ref{pp} together with Lemma~\ref{sets}  (i),
(iii).
\endpf

Let $ \G=\< T_a, T_b \>$ as before. Consider $supp(\G)=S_{a \cup
b}$. If $supp(\G)$ is not a connected surface, and $U$ is one of
its components, we can look at the group $\G|_U$. Certainly if
$\G|_U\cong \mathbb F_2$ then $\G \cong \mathbb F_2$ as well.
Notice that
 an element $g|_U \in \G|_U$ is obtained by dropping the twists in curves which can be isotoped off $U$ from $g \in \G$. So let us characterize the groups $\G$ such that $supp(\G)$ is connected.

\begin{remark}\label{connsupp} {\rm{ Let $\G=\<T_a, T_b\>$ where $T_a,T_b$ are multi-twists. If $supp(\G)$ is connected then $(a_i,b)>0$ and $(a,b_j)>0$ for all $i,j$. }}
\end{remark}

\begin{theorem}\label{2multi2}
For two positive multi-twists
 $T_a=T^{m_1}_{a_1}\cdots T^{m_k}_{a_k}$ and $T_b=T^{n_1}_{b_1}\cdots T^{n_l}_{b_l}$ on the surface $S$, let $\G=\< T_a, T_b\>$  and assume that $supp(\G)$ is connected. Then $\G \cong \mathbb F_2$ except possibly when either
\begin{enumerate}
\item there is $1 \le i \le k$ such that $m_i(a_i,b)=1$ and there is $1 \le j \le l$ such that $n_j(a,b_j)\le 3$, or
\item there is $1 \le j \le l$ such that $n_j(a,b_j)=1$ and there is $1 \le i \le k$ such that $m_i(a_i,b)\le 3$.
\end{enumerate}
\end{theorem}

\proof Suppose that neither of the two cases happen. The group $\G\cong
\mathbb F_2$ if $m_i(a_i,b) \ge 2 $ and $n_j(a,b_j) \ge 2$ for all
$i,j$, by Theorem~\ref{2multifree}. To understand the other cases,
without loss of generality assume that $m_1(a_1,b)=1$. By
Remark~\ref{connsupp} $(a_i, b)>0$ for all $i$ and $(a,b_j)>0$ for
all $j$ since $supp(\G)$ is connected.

 Now put $\l=2$ in Lemma~\ref{sets}. Clearly the condition $m_i(a_i,b) \ge 2 \l\i=1$ is satisfied, so if $n_j(a,b_j) \ge 2\l=4$ for all $j$, using Lemma~\ref{pp} we get $\G \cong \mathbb F_2$. \endpf

One can completely answer the question ``when is a group generated
by powers of Dehn twists isomorphic to $\mathbb{F}_2$?", as
follows:

 \begin{theorem}\label{2twistsfree} Let $A=\{a,b\}$ be a set of two
 simple closed curves on a surface $S$ and $m,n>0$.
Put $\G=\langle T_a^{m}, T_b^{n}\rangle$. The following conditions
are equivalent:
\begin{enumerate}
\item $\G \cong \mathbb F_2$.
\item Either $(a,b)\ge 2$, or
$(a,b)=1$  and
$$\{m,n\} \notin \{ \{1\}, \{1,2\}, \{1,3\} \}.$$
\end{enumerate}
\end{theorem}

\proof By Theorem~\ref{2multi2}, (ii) implies (i). To prove
(i) implies (ii), we must show that for $(a,b)=1$, the groups
$ \langle T_a, T_b^n \rangle$
 are not free for $n=1,2,3$.

\begin{figure}[ht!]\small
$$
  \setlength{\unitlength}{0.045in}
  \begin{picture}(0,0)(0,11)
  %\multiput(0,0)(5,0){25}{\line(0,1){150}}
   %\multiput(0,0)(0,5){25}{\line(1,0){150}}
  \put(42,27){$\delta$}
  \put(18,8){$b$}
  \put(3,25){$a$}
   \end{picture}
   \includegraphics[width=1.8in]{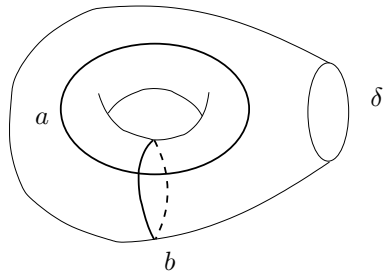}
$$
 \caption{$(T_aT_b)^6=T_\delta$}
     \label{s11}
     \end{figure}

Let us denote $T_a$ by $a$ and $T_b$ by $b$ for brevity. We know
that $(ab)^6$ commutes with both $a$ and $b$, (see
Figure~\ref{s11}; for a proof of this relation see \cite{I2}.) so
 the case
$n=1$ is non-free. Also, notice the famous braid relation
$aba=bab$ (see, for instance~\cite{I2}). Now consider the case
$n=2$. Observe that
$$(ab^2)^2=ab^2 a b^2= ab(bab)b=ab(aba)b=(ab)^3,$$
so $(ab^2)^4  =(ab)^6$ is in the center of $\<a,b^2\>$. In the
case $n=3$, notice that
\begin{eqnarray}
(ab^3)^3&=&ab^3ab^3ab^3  \nonumber \\
&=&ab^2(bab)b(bab)b^2 \nonumber \\
&=&ab^2abababab^2  \nonumber \\
&=&ab(bab)(aba)(bab)b  \nonumber \\
&=&ab(aba)(bab)(aba)b \nonumber \\
&=&(ab)^6. \nonumber
\end{eqnarray}
Therefore $(ab^3)^3$ is in the center of $\<a,b^3\>$.
\endpf

\begin{remark} {\rm After the completion of this work the author
 learned that the isomorphism $\<T_a,T_b\> \cong \mathbb F_2$ for
 $(a,b) \ge 2$ was proved earlier by Ishida  \cite{ISH}.   }
\end{remark}

Let $T_a$, $T_b$ be two positive multi-twists. In the rest of this
section we answer the question ``which words in $\langle T_a,
T_b \rangle $ are relatively pseudo-Anosov?''(see
Definition~\ref{relativepa}).

An element  $f \in \calm(S)$ is called {\it{pure}}~\cite{I} if for
any simple closed curve $c$, $f^n(c)=c$ implies $f(c)=c$. In other
words, by Thurston classification~\cite{FLP}, there is a finite
(possibly empty) set $C=\{c_1,\cdots,c_k\}$ of disjoint simple
closed curves $c_i$ such that $f(c_i)=c_i$ and $f$ keeps all
components of $S \setminus (c_1 \cup \cdots \cup c_k) $ invariant,
and is either identity or pseudo-Anosov on each such component. A
subgroup of $\calm(S)$ is called pure if all elements of it are
pure. Ivanov showed that $\calm(S)$ contains finite-index pure
subgroups, namely, $ker(\calm(S) \to H_1(S,\mathbb Z / m \mathbb Z))$ for $m
\ge 3$ \cite{I}. A relatively pseudo-Anosov word induces a pure element of the mapping class group.

\begin{theorem}\label{2multipa} For two positive multi-twists
 $T_a=T^{m_1}_{a_1}\cdots T^{m_k}_{a_k}$ and
 $T_b=T^{n_1}_{b_1}\cdots T^{n_l}_{b_l}$ on
  the surface $S$, the group $\G=\< T_a, T_b\>$ is pure and relatively
   pseudo-Anosov if any of the following conditions  holds:

 \begin{enumerate}
 \item  For all $i$, $m_i(a_i,b)\ge 2$ and for all $j$, $n_j(a,b_j) \ge 3$.
\item  For all $i$, $m_i(a_i,b)\ge 3$ and for all $j$, $n_j(a,b_j) \ge 2$.
\item  For all $i$, $m_i(a_i,b)\ge 1$ and for all $j$, $n_j(a,b_j) \ge 5$.
\item  For all $i$, $m_i(a_i,b)\ge 5$ and for all $j$, $n_j(a,b_j) \ge 1$.
\end{enumerate}
\end{theorem}

\proof We use Lemma~\ref{tpp} together with Lemma~\ref{sets}
(ii),(iv). First assume that $\l=1+\e$ where $\e$ is a small
irrational number. Notice that $X=N_{a,\l} \cup N_{b,\l\i}$. If
all $m_i(a_i,b)\ge 2 > 2 \l\i$ and $n_i(a,b_i) \ge 3> 2 \l$, one
can use Tower ping-pong to show that if a simple closed curve
intersects $supp(\G)$ then it cannot be mapped to itself by any
element of $\G$ except conjugates of powers of $T_a$ and $T_b$,
which are already known to be pure and relatively pseudo-Anosov.
This proves (i) (A relatively pseudo-Anosov word induces a pure element). Similarly by using $\l=1-\e$, $\e \in \mathbb R
\setminus \mathbb Q$ in Lemma~\ref{sets} (ii),(iv), we get
(ii). To get parts (iii),(iv) we can set $\l=2+\e$ and
$\l=1/2-\e$ respectively and argue similarly.
\endpf

%
%Most of these pure groups are not contained in
%the ones discovered by Ivanov, i.e., of the form $ker(\calm(S) \to
%H_1(S, t \Bbb Z))$, $t \ge 3$, if at least one of the simple
%closed curves $a,b$ is non-separating.

This in particular proves:

\begin{corollary}[Thurston~\cite{FLP}] If $a,b$ are two simple closed curves, which fill
up the closed surface $S$ of genus $g \ge 2$,
 then $\langle T_a,T_b \rangle \cong  \mathbb{F}_2$
 and all elements not conjugate to the powers of $T_a$ and $T_b$
are pseudo-Anosov.
\end{corollary}

\proof If $a,b$ fill up $S$ we must have $(a,b)\ge 3$. Now we can
use Theorem~\ref{2multipa}. \endpf

\begin{theorem}\label{pa(a,b)=2} Let $A=\{a,b\}$ be a set of two
 simple closed curves on a surface $S$ and $m,n>0$ be integers and
  $\G=\langle T_a^{m}, T_b^{n}\rangle$. The following conditions
are equivalent:
\begin{enumerate} \item $\G$ is relatively
pseudo-Anosov.
\item  Either $(a,b)\ge 3$, or $(a,b)=2$ and
$(m,n)\ne (1,1)$, or $(a,b)=1$  and
$$\{m,n\} \notin \{ \{1\}, \{1,2\}, \{1,3\}, \{1,4\}, \{2\}  \}. \ $$
\end{enumerate}
\end{theorem}

\proof If $(a,b) \ge 3$, then $\G$ is relatively pseudo-Anosov for
all $m,n>0$ by Theorem~\ref{2multipa}. If $(a,b)=2$ then $\G$ is
relatively pseudo-Anosov if $m>1$ or $n>1$ by
Theorem~\ref{2multipa}. We prove that if $(a,b)=2$ then
$\G=\<T_a,T_b\>$ is not relatively pseudo-Anosov. We consider two
cases.

\rk{Case 1} $(a,b)=2$ and the algebraic
intersection number of $a,b$ is $\pm 2$.

In this case both $a,b$
can be embedded in a twice punctured torus subsurface of $S$ (see
Figure~\ref{abins120}). We will prove in Proposition~\ref{tbta2}
that $(T_bT_a)^2$ is in fact a multi-twist.

\setlength{\unitlength}{0.045in}

 \begin{figure}[ht!]\small
 $$\begin{picture}(0,0)(0,11)
   %\multiput(0,0)(5,0){25}{\line(0,1){150}}
   %\multiput(0,0)(0,5){25}{\line(1,0){150}}
 \put(-5,2){$a$}
  \put(0,-24){$b$}
  \put(-21,-1){$\partial_2$}
  \put(19.5,-1){$\partial_1$}
   \end{picture}$$

  {\centerline{ \includegraphics[width= 2.7 in]{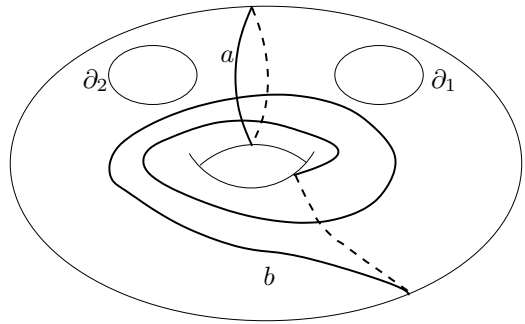}}}
     \caption{ $S_{1,2,0}$ and the curves $a$ and $b$}
     \label{abins120}
     \end{figure}

\rk{Case 2} $(a,b)=2$ but the algebraic
intersection number of $a,b$ is $0$. 

In this case $a,b$ can be embedded in a 4-punctured sphere. According
to the lantern relation~\cite{I2} (see Figure~\ref{s040}), $T_bT_a$ is
a multi-twist.

\setlength{\unitlength}{0.045in}

 \begin{figure}[ht!]\small
 $$\begin{picture}(0,0)(0,11)
   %\multiput(0,0)(5,0){25}{\line(0,1){150}}
   %\multiput(0,0)(0,5){25}{\line(1,0){150}}
   \put(-16,-7){$a$}
   \put(16,-8){$c$}
 \put(-0.7,-.3){$\partial_2$}
 \put(-15,-25){$\partial_3$}
 \put(12,-25){$\partial_4$}
  \put(0,-36){$b$}
  \put(1,10){$\partial_1$}
   \end{picture}$$

  {\centerline{ \includegraphics[width=2.25 in]{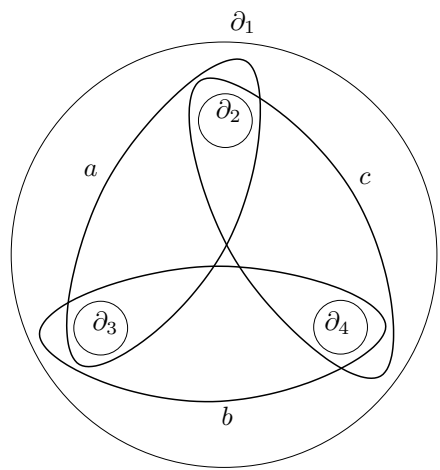}}}
     \caption{$T_aT_bT_c=T_{\partial_1}T_{\partial_2}
T_{\partial_3}T_{\partial_4} $}
     \label{s040}
     \end{figure}

This proves that when $(a,b)=2$, $\<T_a,T_b \>$ is not relatively
pseudo-Anosov.

 If $(a,b)=1$, the group $\G$ is relatively pseudo-Anosov except
 possibly when $\{m,n\}$ is one of $\{1,i\}$, $i=1,2,3,4$, or
 $\{m,n\}=\{2,2\}$, by Theorem~\ref{2multipa}.
The groups $\<T_a,T_b\>$, $\<T_a, T_b^2\>$ and $\<T_a, T_b^3\>$
are not relatively  pseudo-Anosov because the map
$(T_aT_b)^6=(T_aT_b^2)^4=(T_aT_b^3)^3$ is in fact a Dehn twist in
the boundary of the surface defined by $a,b$
 (see Figure~\ref{s11}) hence they induce the identity on $S_{a,b}$.

If $(a,b)=1$, then $\G=\< T_a^2, T_b^2 \>$ is not relatively
pseudo-Anosov. This is because $T_b^2T_a^2$ has a trace of $-2$
and hence is reducible (see Remark~\ref{mcg1}). Similarly when $(a,b)=1$,
the maps $T_aT_b^4$ and $T_a^4T_b$ both have a trace of $-2$ and hence are reducible.
\endpf

We saw that if $a,b$ are two simple closed curves with $(a,b) \ge
2$,  a word $w(T_a,T_b) \in \<T_a,T_b \>$ is relatively
pseudo-Anosov except possibly when $(a,b)=2$. In the following
theorem we narrow down the search for words which are not relatively pseudo-Anosov in this case.

\begin{theorem}\label{(a,b)=2} Let $a,b$ be two simple closed curves on a surface
$S$ with $(a,b)=2$. Then a word $w$ in $T_a, T_b$ representing an
element of $ \<T_a, T_b\>$ is a pure and relatively pseudo-Anosov
unless possibly when $w$ is cyclically reducible to a power of
either  $T_bT_a^{-1}$ or $T_bT_a$.
\end{theorem}

\proof The proof is based on repeated application of
Lemma~\ref{flplemma}. Clearly $T_a$ and $T_b$ are both pure and
relatively pseudo-Anosov. So in what follows we assume that $w$ is
a cyclically reduced word of length $>1$. Let $X,A, N_{a,1},
N_{b,1}$ be defined as in the beginning of this section. Let
$$Y=\{ x \in X \ | \ (x,a)=(x,b) \}.$$
Hence $X$ is a disjoint union $N_{a,1} \cup N_{b,1} \cup Y$.

 By Lemma~\ref{sets}, we have $T^{\pm
n}_a(N_{b,1}) \subset N_{a,1}$ and
 $T^{\pm n}_b(N_{a,1}) \subset N_{b,1}$ for all $n>0$. Moreover, for $x \in N_{b,1}$,
we have $||T_a^{\pm n}(x)|| > ||x||$ and for $x \in N_{a,1}$, $||T_b^{\pm n}(x)||>||x||$.
 By the same lemma, $T^{\pm
n}_a(Y) \subset N_{a,1}$ and
 $T^{\pm n}_b(Y) \subset N_{b,1}$ for all $n \ge 2$ (This follows by applying the lemma to $\l=1+\e$
 and $\l=1-\e$, where $\e$ is a small positive number).

 Let $w=T_b^{n_k}T_a^{m_k} \cdots
 T_{b}^{n_1}T_a^{m_1}$ be a cyclically reduced word, where $m_i,n_i \ne 0$ and $k \ge 1$.
  If  any of $m_i$ is greater that
 $1$ in absolute value, we can assume without loss of generality that $|m_1|>1$, by conjugation.
 Therefore if $x \in Y \cup N_{b,1}$,  then $T_a^{m_1}(x) \in
 N_{a,1}$ and hence $||w^n(x)||>||x||$ for all $n>0$. Hence $w^n(x)
 \ne x$ for all integers $n$. If $x  \in N_{a,1}$, then
 $||w^{-n}(x)||>||x||$ for $n>0$, and so $w^n(x) \ne x$ for all
 $n$. This shows that $w$ is relatively
 pseudo-Anosov. The case where some $|n_i|>1$ follows by symmetry
 by replacing $w$ with $w\i$.

 So let us assume that for all $1 \le i \le  k$, we have
 $m_i,n_i=\pm 1$. If $w$ is not conjugate to a power of $T_bT_a$
 or $T_bT_a^{-1}$, by conjugating $w$ we can assume either
$m_1 \ne m_2$, or
 $n_k \ne n_{k-1}.$
We assume the former. The latter can be dealt with similarly by
symmetry and replacing $w$ with $w\i$. In this case the word $w$
could have any of the following forms:

\begin{enumerate}
  \item   $w=T_b^{n_k}T_a^{m_k} \cdots T_a T_b T_a\i$,
  \item   $w=T_b^{n_k}T_a^{m_k} \cdots T_a\i T_b T_a$,
  \item   $w=T_b^{n_k}T_a^{m_k} \cdots T_a T_b\i T_a\i$,
  \item   $w=T_b^{n_k}T_a^{m_k} \cdots T_a\i T_b\i T_a$.
\end{enumerate}

Suppose, for example, that  $w=T_b^{n_k}T_a^{m_k} \cdots T_a T_b
T_a\i$. As before, if $x \in N_{a,1} \cup N_{b,1}$, we get that
$w^n(x) \ne x$ for all $n>0$. So let us assume that $x \in Y$.
Then, by definition of $Y$, $(x,a)=(x,b)=p>0$. Then we have
$(T_a\i(x),a)=p$ and by Lemma~\ref{flplemma},
$$|(T_a\i(x),b)-(a,b)(x,a)|\le (x,b),$$
which implies  $p \le (T_a^{-1}(x),b) \le 3p$. If $p
<(T_a^{-1}(x),b)$, then $T_a\i(x) \in N_{a,1}$ and so $w^n(x) \ne
x$, for all $n>0$. So let us assume $(T_a^{-1}(x),b)=p$.  Notice
that this implies $(T_b T_a\i(x), b)=p$. Again by
Lemma~\ref{flplemma},
$$|(T_b( T_a\i(x)), a)-(a,b)(b, T_a\i(x))| \le (T_a\i(x),a),$$
which gives $p \le (T_b T_a\i(x), a)  \le 3p$. Again, if $p < (T_b
T_a\i(x),a)$, then $T_bT_a\i(x) \in N_{b,1}$ which implies $w^n(x)
\ne x$ for $n >0$.  Otherwise, we can further assume that $(T_b
T_a\i(x),a)=p$. Notice that this gives $(T_aT_b T_a\i(x),a)=p$. At
this point it looks like the argument is going to go on forever,
but here is a new ingredient. For any mapping class $f$, we have
the following well-known equation: $fT_bf\i=T_{f(b)}$. In
particular: $T_aT_bT_a\i=T_{T_a(b)}$.

\medskip
$\bullet$\qua{\bf Claim}\qua $(T_a(b),b)=4$
\medskip

This follows from Lemma~\ref{flplemma}:
$|(T_a(b),b)-(b,a)(a,b)|\le (b,b)$.

Now by the same lemma,
$$|( T_{T_a(b)}(x),b)-(T_a(b),x)(T_a(b),b)| \le (x,b),$$
which gives $ 3p \le (T_aT_b T_a\i(x),b) \le 5p$, i.e., $T_aT_b
T_a\i(x) \in N_{a,1}$, and so $w^n(x) \ne x$ for all $n \ge 0$.
The other cases (ii),(iii) and (iv) follow similarly. \endpf

\section{The case of two simple closed curves with intersection
number 2 filling a 4-punctured sphere}\label{S20}

Let $a,b$ be two simple closed curves such that  $(a,b)=2$ and
$S_{\{a,b\}}$ is a four-holed sphere. (Figure~\ref{s040}).

The relation $T_aT_bT_c=T_{\partial_1}T_{\partial_2}
T_{\partial_3}T_{\partial_4} $ was discovered by Dehn~\cite{D} and
later on by Johnson~\cite{J}. A proof of the lantern relation can
be found in~\cite{I2}. Note the commutativity between the various
twists.

\begin{prop}\label{paforcase1} In the group $\<T_a,T_b\>$ all words are pure. All
words are relatively pseudo-Anosov except precisely words that are
cyclically reducible to a non-zero power of $T_bT_a$.
\end{prop}

\proof The lantern relation implies:
$$T_aT_b=T_c^{ -1}T_{\partial_1}T_{\partial_2}
T_{\partial_3}T_{\partial_4}.$$ This shows that $T_aT_b$ (and
hence its conjugate $T_bT_a$) is a multi-twist. Notice that
$$T_a\i T_b=T_a^{-2}T_c\i T_{\partial_1}T_{\partial_2}
T_{\partial_3}T_{\partial_4}.$$ Hence restricted to $S_{a,b}$,
$T_a\i T_b=T_a^{-2}T_c\i$. But the group $\< T_a^2, T_c \>$ is
pure and relatively pseudo-Anosov by Theorem~\ref{pa(a,b)=2},
which shows that $T_a\i T_b$ (and hence its conjugate $T_bT_a\i$)
is pure relatively pseudo-Anosov. Moreover, $a,c$ fill the same
surface as $a,b$. Finally, we invoke Theorem~\ref{(a,b)=2}.
\endpf

\section{The case of two simple closed curves with intersection
number 2 filling a twice-punctured torus}\label{S12}

Let $S_{g,b,n}$ denote a surface of  genus $g$ with $b$ boundary
components and $n$ punctures. Let $a$ and $b$ be two simple
closed curves such that $(a,b)=2$ and assume both intersections
have the same sign. In this case $a$ and $b$ are both
non-separating. One can therefore assume, up to diffeomorphism that they
are as given in Figure~\ref{abins120}. Since the regular
neighborhood of $a \cup b$ is homeomorphic to $S_{1,2,0}$, the
surface filled by $a,b$ is $S_{1,i,j}$, for $i,j=0,1,2$, $i+j\le 2$.

\setlength{\unitlength}{0.045in}

 \begin{figure}[ht!]\small
 $$\begin{picture}(0,0)(0,11)
   %\multiput(0,0)(5,0){25}{\line(0,1){150}}
   %\multiput(0,0)(0,5){25}{\line(1,0){150}}
   \put(-13,-7){$\gamma$}
 \put(-0.7,-.3){$T_a$}
 \put(-5,-11){$T_b$}
 \put(-.7,-26){$T_a$}
 \put(-5,-37){$T_b$}
   \put(-13,-32){$\gamma\p$}
   \put(-13,-59){$\gamma$}
   \put(-26.5,6){$a$}
   \put(26,-11){$b$}
   \end{picture}$$

  {\centerline{ \includegraphics[width= 3.6 in]{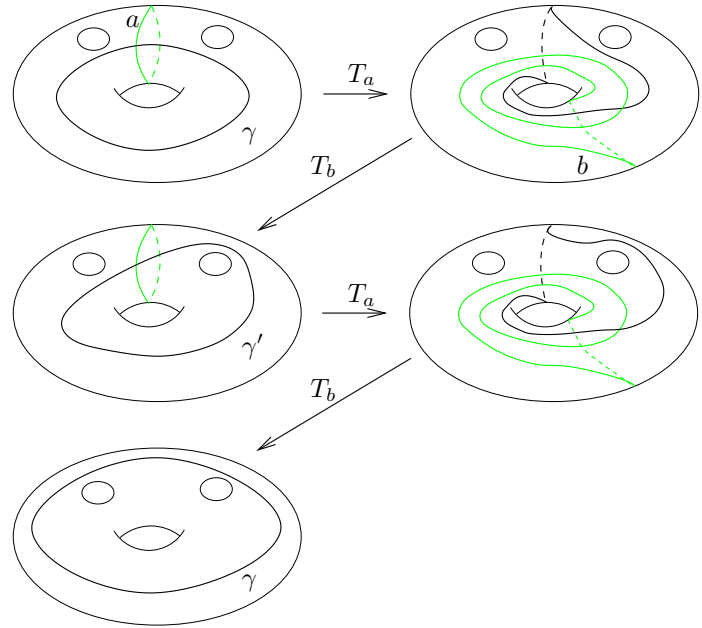}}}
     \caption{$(T_bT_a)^2(\gamma)=\gamma$ and $(T_bT_a)(\gamma)=\gamma\p$}
     \label{gamma}
     \end{figure}

Assume that $S_{\{a , b\}}=S_{1,2,0}$. Let $\gamma$ and $\gamma\p$
be the curves defined in Figure~\ref{gamma}. By following
Figure~\ref{gamma} one can see that
$(T_bT_a)^2(\gamma)=\gamma$, preserving the orientation. Since by
definition $T_bT_a(\gamma)=\gamma\p$, one also gets
$(T_bT_a)^2(\gamma\p)=\gamma\p$. Now notice that $\gamma$ and
$\gamma\p$ cut up $S_{1,2,0}$ into two pairs of pants. Hence
$(T_bT_a)^2$ is a multi-twist in curves $\gamma, \gamma\p$,
$\partial_1$ and $\partial_2$.

\begin{prop}\label{tbta2}
With the notation in Figures~\ref{abins120} and \ref{gamma}, we
have $$(T_bT_a)^2=T_{\partial_1}T_{\partial_2}T_{\gamma}^{-4}
T_{\gamma\p}^{-4}.$$
\end{prop}

\setlength{\unitlength}{0.045in}

 \begin{figure}[ht!]\small
 $$\begin{picture}(0,0)(0,11)
   %\multiput(0,0)(5,0){25}{\line(0,1){150}}
   %\multiput(0,0)(0,5){25}{\line(1,0){150}}
\put(-18,-1){$I$} \put( -2, 0){$T_a$} \put( -2, -25.5){$T_a$}
\put( -2, -51){$=$} \put( -7,- 11.5){$T_b$} \put( -7,-
37.5){$T_b$}
   \end{picture}$$

  {\centerline{ \includegraphics[width= 3.6 in]{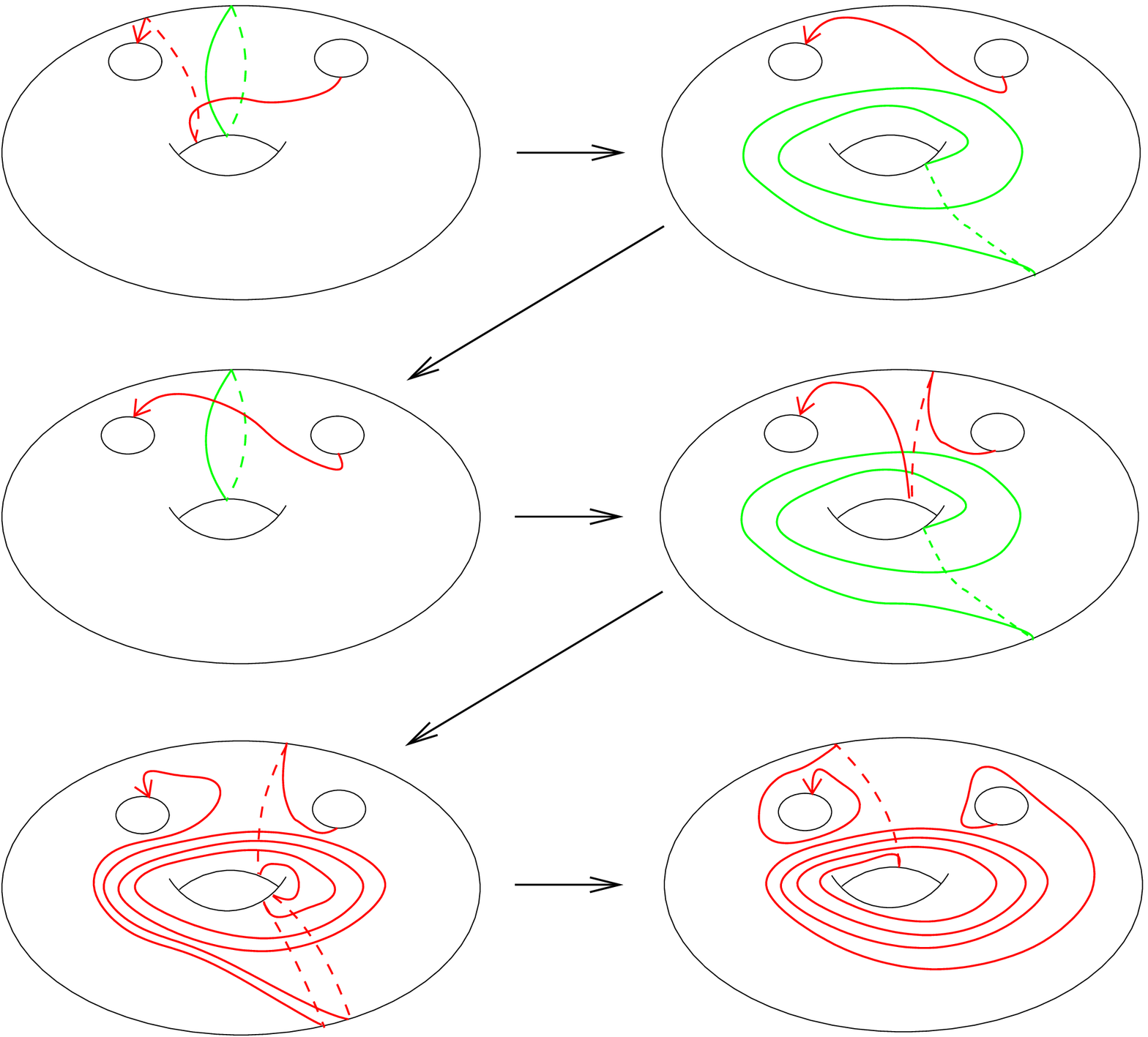}}}
     \caption{The arc $I$, and $(T_bT_a)^2(I)$}
     \label{arcI}
     \end{figure}

\proof Since $(T_bT_a)^2$ fixes $\partial_1, \partial_2, \gamma$ and $\gamma\p$, it has to be a multi-twist
 in these curves. We consider an arc joining $\partial_1$ to $\partial_2$ crossing $\gamma$
 once as in Figure~\ref{arcI}. We apply $(T_bT_a)^2$ to $I$, and the result is the same as applying
 $T_{\partial_1}T_{\partial_2}T_{\gamma}^{-4}$ to $I$ (again see Figure~\ref{arcI}). Hence
$$(T_bT_a)^2=T_{\partial_1}T_{\partial_2}T_{\gamma}^{-4}T_{\gamma\p}^n,$$
where $n$ is to be found. One can argue by drawing another arc joining $\partial_1$ to $\partial_2$
passing through $\gamma\p$ once, but here is a simpler way: We know that $(T_bT_a)(\gamma)=\gamma\p$ and
 $(T_bT_a)(\gamma\p)=\gamma$
 (see Figure~\ref{gamma}), so if we conjugate the above equation by $T_bT_a$, we get:
$$(T_bT_a)^2=T_{\partial_1}T_{\partial_2}T_{\gamma\p}^{-4}T_{\gamma}^n,$$
which shows that $n=-4$. \endpf

\begin{prop}\label{TbTa\i-S120} The word $T_bT_a\i$ is relatively pseudo-Anosov.
\end{prop}

\proof We use a ``brute force'' method to show that restricted to
$S=S_{1,0,2}$ (the boundary components $\partial_1, \partial_2$
shrunk to punctures $p_1,p_2$, respectively), the word $T_bT_a\i$
induces a pseudo-Anosov map. It is enough to show that the word
$T_b\i T_a$ is relatively pseudo-Anosov. Let $f$ be the mapping
class induced by the word $T_b\i T_a$. We will find measured
laminations $\calf_1 , \calf_2$ and $\l
>0$ such that $f(\calf_1) =\l \calf_1$ and $f(\calf_2) = \l\i
\calf_2$. To this end, we use the theory of measured train-tracks.
For a review of these methods and the theory, see for
example~\cite{ht}.

\setlength{\unitlength}{0.045in}

 \begin{figure}[ht!]\small
 $$\begin{picture}(0,0)(0,11)
   %\multiput(0,0)(5,0){25}{\line(0,1){150}}
   %\multiput(0,0)(0,5){25}{\line(1,0){150}}
\put(28,-1){$a$} \put( -2, -27){$b$} \put( -7,- 10){$b$} \put(
-25,- 37.5){$b$} \put( 39,11){$p_1$} \put( 0, 11){$p_2$}
\put(-39,11){$p_1$} \put(0,-38){$p_2$} \put( 41,-40){$p_1$} \put(-
42,-40){$p_1$}
   \end{picture}$$

  {\centerline{ \includegraphics[width= 3.6 in]{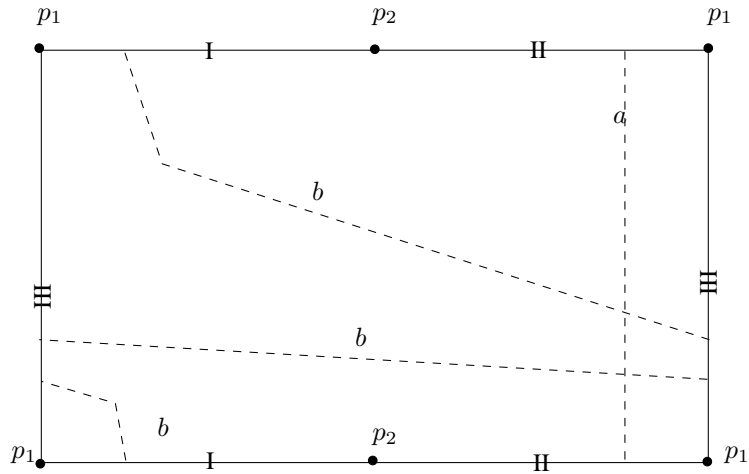}}}
     \caption{The surface $S_{1,0,2}$ cut-open into a polygon $R$}
     \label{polygon}
     \end{figure}

\setlength{\unitlength}{0.045in}

 \begin{figure}[ht!]\small
 $$\begin{picture}(0,0)(0,11)
   %\multiput(0,0)(5,0){25}{\line(0,1){150}}
   %\multiput(0,0)(0,5){25}{\line(1,0){150}}
\put( 15, -32){$x$} \put( 15,0){$x$} \put( -20,0){$y$} \put(
-20,-32){$y$} \put( -3,-15){$z$}

\put( 39,11){$p_1$} \put( 0, 11){$p_2$} \put(-39,11){$p_1$}
\put(0,-38){$p_2$} \put( 41,-40){$p_1$} \put(- 42,-40){$p_1$}
   \end{picture}$$

  {\centerline{ \includegraphics[width= 3.6 in]{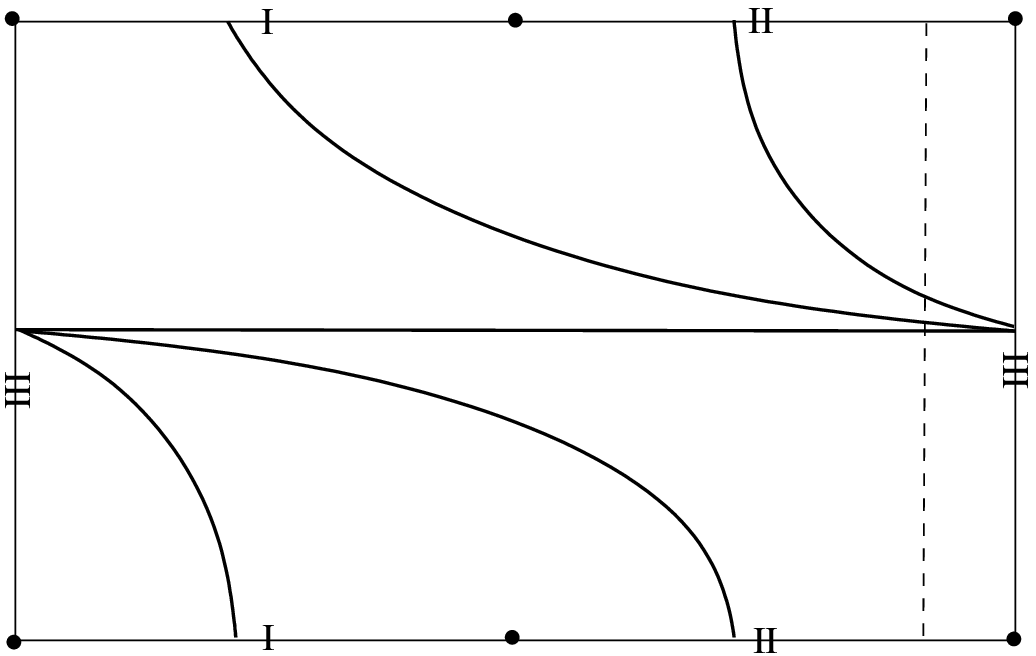}}}
     \caption{The measured train-track $\tau(x,y,z)$}
     \label{tau}
     \end{figure}

\setlength{\unitlength}{0.045in}

 \begin{figure}[ht!]\small
 $$\begin{picture}(0,0)(0,11)
   %\multiput(0,0)(5,0){25}{\line(0,1){150}}
   %\multiput(0,0)(0,5){25}{\line(1,0){150}}
\put( 14,-7){$x$} \put( 23,-7){$x+y+z$} \put( -18,0){$y$} \put(
-22,-32){$y$} \put( -11,-19){$x+z$}

\put( 39,11){$p_1$} \put( 0, 11){$p_2$} \put(-39,11){$p_1$}
\put(0,-38){$p_2$} \put( 41,-40){$p_1$} \put(- 42,-40){$p_1$}
   \end{picture}$$

  {\centerline{ \includegraphics[width= 3.6 in]{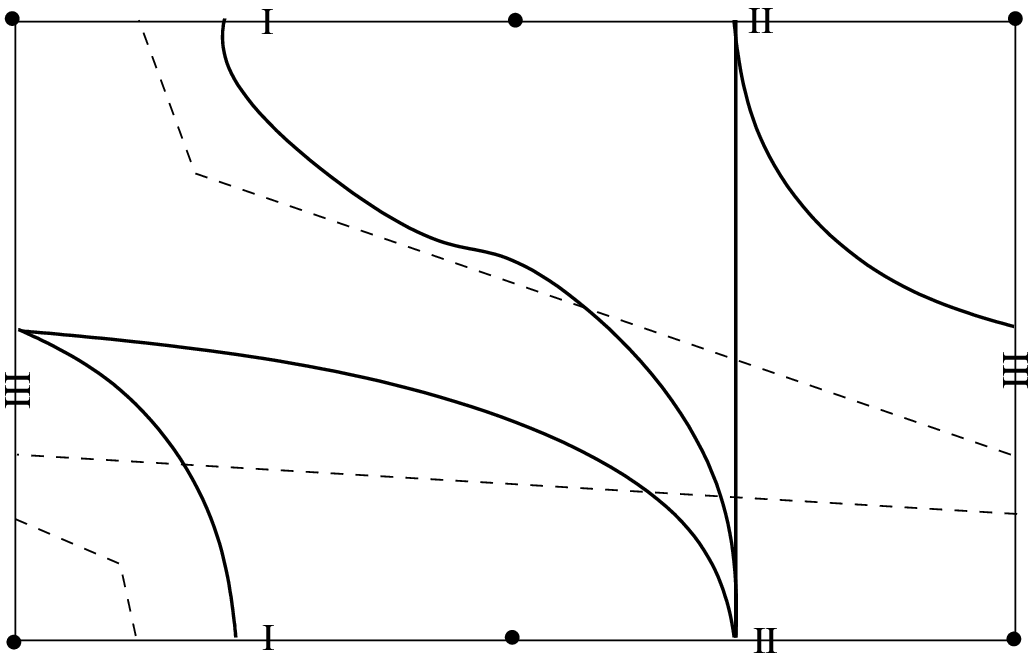}}}
     \caption{The measured train-track $T_a(\tau(x,y,z))$}
     \label{tau1}
     \end{figure}

\setlength{\unitlength}{0.045in}

 \begin{figure}[ht!]\small
 $$\begin{picture}(0,0)(0,11)
   %\multiput(0,0)(5,0){25}{\line(0,1){150}}
   %\multiput(0,0)(0,5){25}{\line(1,0){150}}
\put( 23,-5){$2x+y+z$} \put( -31,-3){$3x+4y+z$} \put(
-11,-15){$3x+3y+z$} \put( -22,-33){$3x+4y+z$} \put( 39,11){$p_1$}
\put( 0, 11){$p_2$} \put(-39,11){$p_1$} \put(0,-38){$p_2$} \put(
41,-40){$p_1$} \put(- 42,-40){$p_1$} \put( 16,-33){$2x+y+z$}
   \end{picture}$$

  {\centerline{ \includegraphics[width= 3.6 in]{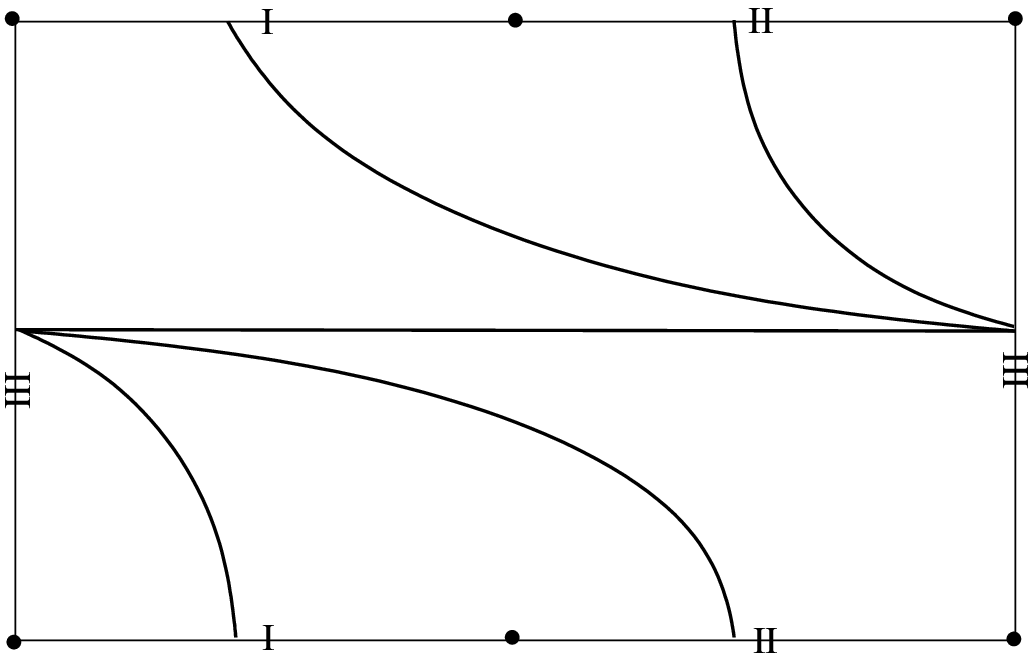}}}
     \caption{The measured train-track $f(\tau(x,y,z))={T_b}\i {T_a}(\tau(x,y,z))$}
     \label{tau2}
     \end{figure}

Consider the polygon $R$ obtained by cutting $S$ open as in
 Figure~\ref{polygon}.
 Identifying parallel sides of  $R$ yields back the surface $S$. Consider
the measured train-track $\tau=\tau(x,y,z)$ ($x,y,z \ge 0$) on $S$
defined as in Figure~\ref{tau}. We can calculate the image
$f(\tau)$ in two steps as in Figures~\ref{tau1} and~\ref{tau2}.
(Remember that Dehn twists are right-handed.) Luckily the action
of $f$ on the space of measures on $\tau$ is linear, so we can
easily find fixed laminations carried on $\tau$: The matrix
representing $f$ on the space of measured laminations carried on
$\tau$  is
 $$\begin{pmatrix} 2 & 3 & 3 \\  1 & 4 & 3 \\  1 & 1 & 1  \end{pmatrix}. $$
  This matrix has eigenvalues $1,
3+\sqrt{10}, 3-\sqrt{10}$. The only eigenvalue that has a
non-negative eigenvector is $\l=3+\sqrt{10}$ and the eigenvector
corresponds to the measured train-track $\tau(2+\sqrt{10},
2+\sqrt{10},2)$ (up to a positive factor). If we ``fatten up" this
measured train-track, we get a lamination $\calf_1$ as in
Figure~\ref{fat1}
 with the property $f(\calf_1)=\l \calf_1$. Notice that,
geometrically, all leaves have slope -1. One can see that there
are no closed loops of leaves (if there were they would have been
caught as eigenvectors already). Also, there is no leaf in
$\calf_1$ connecting a puncture to a puncture, since if it were so
$\sqrt{10}$ would be rational. Similarly, one can find a
lamination $\calf_2$ which satisfies $f(\calf_2) =\l\i \calf_2$.
But establishing such $\calf_1$ is already enough to show that $f$
is pseudo-Anosov on $S$, hence proving the proposition. \endpf

\setlength{\unitlength}{0.045in}

 \begin{figure}[ht!]\small
 $$\begin{picture}(0,0)(0,11)
   %\multiput(0,0)(5,0){25}{\line(0,1){150}}
   %\multiput(0,0)(0,5){25}{\line(1,0){150}}

\put( 25,11){$p_1$} \put( 0, 11){$p_2$} \put(-25,11){$p_1$}
\put(-27,5){$2$} \put(26.5,-39){$2$} \put(0,-40){$p_2$}
\put(27,-42){$p_1$} \put(-27,-42){$p_1$}
 \put(7,11){$2+\sqrt{10}$} \put( -16,11){$2+\sqrt{10}$}
\put(27,-2){$2+\sqrt{10}$} \put(27,-24){$2+\sqrt{10}$}
\put(-37,-10){$2+\sqrt{10}$} \put(-37,-32){$2+\sqrt{10}$}
 \put(5.5,-40){$2+\sqrt{10}$}
\put(-17,-40){$2+\sqrt{10}$}  \end{picture}$$

  {\centerline{ \includegraphics[width= 2.25 in]{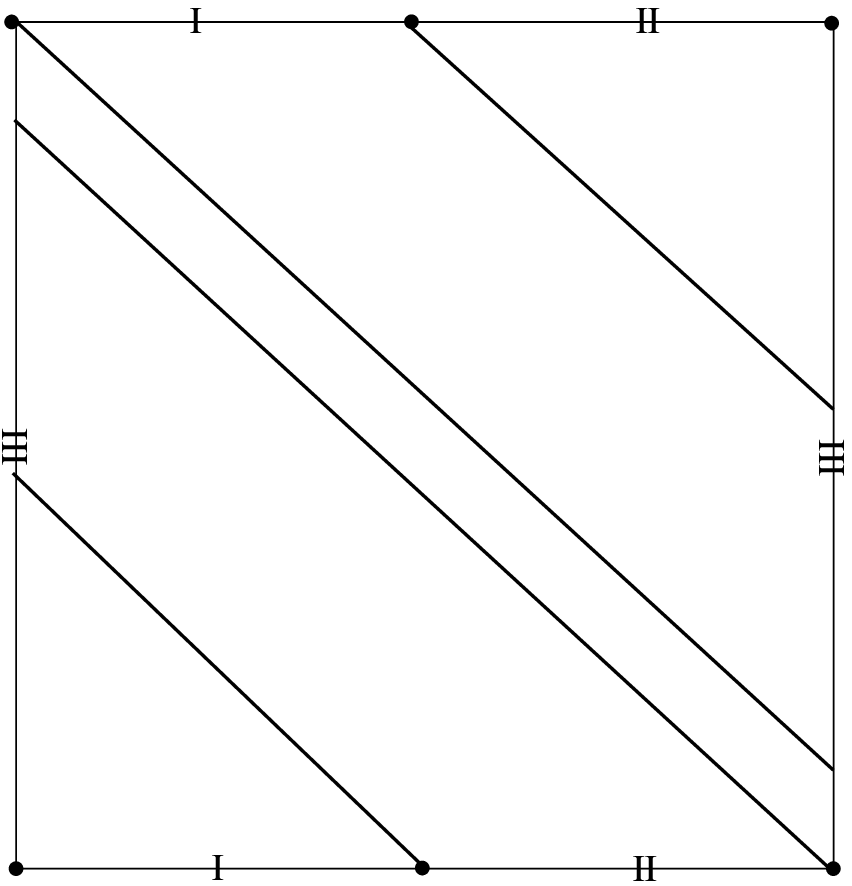}}}
     \caption{The measured lamination $\calf_1$ satisfies $f(\calf_1)= \l \calf_1$}
     \label{fat1}
     \end{figure}

\begin{corollary}\label{paforcase2} All words in $\<T_a, T_b\>$ are pure except precisely the ones
 conjugate to the odd powers of $T_bT_a$. All words are relatively pseudo-Anosov except
precisely the ones cyclically reducible to $(T_bT_a)^n$ for some
non-zero integer $n$.
\end{corollary}

\proof  Notice that the words conjugate to powers of $T_bT_a\i$ are
relatively pseudo-Anosov and hence pure. Now the claim follows
from Theorem~\ref{(a,b)=2} and Proposition~\ref{TbTa\i-S120}.
\endpf

\section{Application to Lantern-type relations}\label{applan}

\begin{theorem}\label{lantern} Let $a,b$ be two simple closed curves on a surface
$S$ such that $(a,b) \ge 2$. Let $w$ be a word in $T_a,T_b$ which
is not cyclically reducible to a power of $T_a$ or $T_b$, but
representing an element in $\calm(S)$ which is a multi-twist. Then
$(a,b)=2$ and exactly one of the following conditions hold:

\begin{enumerate}
  \item The curves $a,b$ have algebraic intersection number $0$,
  the word $w$ can be cyclically reduced to $(T_aT_b)^n$
  for some $n \in \mathbb{Z}$, and $$T_aT_b=T_{\partial_1}T_{\partial_2}
T_{\partial_3}T_{\partial_4}T_c^{-1}.$$ (See Figure~\ref{s040}).
  \item  The curves $a,b$ have algebraic intersection number $2$,
  the word $w$ can be cyclically reduced to $(T_bT_a)^{2n}$
  for some $n \in \mathbb{Z}$, and $$(T_bT_a)^2=T_{\partial_1}T_{\partial_2}
T_{\g}^{-4}T_{\g\p}^{-4}.$$ (See
Figures~\ref{abins120},~\ref{gamma}).
\end{enumerate}
\end{theorem}

\proof By Theorem~\ref{2multipa} $(a,b)=2$, because otherwise $w$
will be relatively pseudo-Anosov, with a support which is not a union of annuli, so it cannot be a multi-twist. Now apply
Proposition~\ref{paforcase1} and Corollary~\ref{paforcase2}.
\endpf

\begin{remark}\label{mcg1} \rm{
It is well known that $\calm_{1,0,1}=\textrm{SL}(2,
\mathbb{Z})$. Also,
$$\sltz =\< t=\begin{pmatrix} 1 & 1 \\  0 & 1 \end{pmatrix},
q=\begin{pmatrix} 0 & 1 \\  -1 & 0 \end{pmatrix}\>.$$
We have a short exact sequence~\cite{bir}
 $$0 \to\mathbb  Z \to \calm_{1,1,0} \to \calm_{1,0,1}=
\sltz \to 0.$$
 The Dehn twists $T_a$ and $T_b$ in
Figure~\ref{s11} induce the matrices
$$  t=\begin{pmatrix} 1 & 1 \\  0 & 1 \end{pmatrix} \text{ and }
s= \begin{pmatrix} 1 & 0 \\  -1 & 1 \end{pmatrix} $$
 in $\sltz$,
respectively. Clearly $sts=q$ and hence $\sltz=\<t,s\>$. In this
case, a word in $\<T_a, T_b\>$ is pseudo-Anosov if and only if the
trace of the corresponding matrix has absolute value of $2$ or
more. Such a word is a multi-twist if the corresponding matrix has
trace $2$. }
\end{remark}

\begin{definition} A relation $w(T_a,T_b)=T_C$ is called
lantern-like if  $T_a,T_b$ are Dehn twists, and $T_C$ is a
multi-twist with $C$ having at least $3$ components.
\end{definition}

\begin{theorem}\label{lanternlike} The only lantern-like relations in any mapping
class group are described in Theorem~\ref{lantern}.
\end{theorem}

\proof We have to only show that, if $(a,b)=1$, they cannot form a
lantern-like relation. But in that case, $a,b$ are supported in a
once-punctured torus, hence $T_C$ can be made of twists in the
boundary and at most one simple closed curve in that torus. \endpf

\section{Groups generated by $n \ge 3$ powers of
twists}\label{nge3}

In this section the phrase ``$i \ne j \ne k$'' means that $i,j,k$
are distinct. Let $a_1,\cdots a_n$ be $n \ge 3$ simple closed
curves on a surface $S$ such
 that $(a_i,a_j)>0$ for $i\ne j$.

Let  $\l_{ijk}>1$ and $\mu_{ij}>0$ (for $i\ne j \ne k$) be real
numbers such that $\mu_{ji}=\mu_{ij}\i$. Put $\l=(\l_{ijk})_{i\ne
j \ne k}$ and $\mu=(\mu_{ij})_{i \ne j}$. Define the set of simple
closed curves
$$N_{a_i}=N_{a_i,\l, \mu}=
\{x \ | \ (x,a_i)<\mu_{ij}(x,a_j), \ {(x,a_k) \over
(x,a_j)}<\l_{ijk} { (a_i,a_k) \over (a_i,a_j)},\ \forall j \ne k
\ne i \},$$
for $i=1,\cdots, n$. Note that $a_i \in N_{a_i}$.

\begin{lemma}\label{3.1} Let $a_1,\cdots,a_n$ be a set of $n\ \ge 3$ simple closed
curves such that $(a_i,a_j)\ne 0$ for $i \ne j$.
\begin{enumerate}

\item The sets $N_{a_i}$,
$i=1,\cdots,n$ are mutually disjoint.

\item
 For $1 \le i\ne j \le  n$, we have $T^{\pm \nu}_{a_i}(N_{a_j})
\subset  N_{a_i}$ for
\begin{eqnarray*} \nu \ge \max &\{&\ { 2 \over \mu_{ij}(a_i,a_j)}, \\
&&  {1 \over \mu_{ik}(a_i,a_k)} + \l_{jik} {(a_j,a_k) \over
(a_i,a_j)(a_i,a_k)},
\\
&&{ \l_{jil} \over \l_{ikl}-1}{(a_j,a_l) \over(a_i,a_l)
 (a_j,a_i)}+{\l_{ikl} \l_{jik} \over \l_{ikl}-1}
{(a_j,a_k)
\over (a_j,a_i)(a_i, a_k)}, \\
&& { 1 \over (\l_{ikj}-1)\mu_{ij}(a_i,a_j)}+{\l_{ikj} \l_{jik}
\over \l_{ikj}-1} {(a_j,a_k)
\over (a_j,a_i)(a_i, a_k)}, \\
&&  {\l_{ijl} \over (\l_{ijl}-1)\mu_{ij}(a_i,a_j)} +{\l_{jil}
\over \l_{ijl}-1}{(a_j,a_l) \over (a_j,a_i)(a_i,a_l)} \}_{k\ne
l\ne i}. \end{eqnarray*}

\end{enumerate}
\end{lemma}

\proof (i) is clear. To prove (ii), consider $x \in N_{a_j}$.
We have
$$(T^{\pm \nu}_{a_i}(x),a_j) \ge \nu(a_i,a_j)(x,a_i)-(x,a_j)
> \mu_{ji}(x,a_i)=\mu_{ji}(T^{\pm \nu}_{a_i}(x), a_i)$$
for $\nu \ge{ 2 \mu_{ji} \over (a_i,a_j)}$. Let $k\ne i,j$. Then
$$(T^{\pm \nu}_{a_i}(x),a_k) \ge \nu(a_i,a_k)(x,a_i)-(x,a_k)
>\mu_{ki} (x,a_i)$$
if
$$ \nu \ge {1 \over \mu_{ik}(a_i,a_k)} + \l_{jik} {(a_j,a_k) \over (a_i,a_j)(a_i,a_k)}.$$
Let $k,l\ne i$. Then
$$(T^{\pm \nu}_{a_i}(x), a_l)/( T^{\pm \nu}_{a_i}(x), a_k)<\l_{ikl}
(a_i,a_l)/(a_i, a_k)$$ if and only if
$$ {   (a_i, a_k)(T_{a_i}^{\pm \nu}(x), a_l) }
< \l_{ikl} {(a_i,a_l)  (T_{a_i}^{\pm \nu}(x), a_k)}.$$ This will
hold if
\begin{equation}\label{eqn1}
{ (a_i, a_k)(\nu(a_i,a_l)(x,a_i)+(x,a_l)) } <
 \l_{ikl} {(a_i,a_l)(  \nu(a_i,a_k)(x,a_i)-(x,a_k)) }.
\end{equation}

The inequality~(\ref{eqn1}) is equivalent to
\begin{equation}\label{eqn2}
\nu(a_i,a_l)(\l_{ikl}-1) > {(x,a_l) \over (x,a_i)}+\l_{ikl}
{(x,a_k)(a_i,a_l) \over (x,a_i)(a_i, a_k)}.
\end{equation}
 (One has $(x,a_i)>0$ since $x \in N_{a_j}$.) Therefore for $l\ne
j$ and $k\ne j$, it is enough to have
$$\nu(a_i,a_l)(\l_{ikl}-1) \ge \l_{jil} {(a_j,a_l) \over (a_j,a_i)}+\l_{ikl}
\l_{jik} {(a_j,a_k)(a_i,a_l) \over (a_j,a_i)(a_i, a_k)};$$ i.e.,
$$\nu\ge{ \l_{jil} \over \l_{ikl}-1}{(a_j,a_l) \over(a_i,a_l)
 (a_j,a_i)}+{\l_{ikl} \l_{jik} \over \l_{ikl}-1}
{(a_j,a_k) \over (a_j,a_i)(a_i, a_k)}.$$ If $l=j$ (and so $ k\ne
j$) then one can replace the inequality~(\ref{eqn2}) with
$$\nu(a_i,a_l)(\l_{ikl}-1) \ge \mu_{ji}+\l_{ikl} {(x,a_k)(a_i,a_l)
\over (x,a_i)(a_i, a_k)}$$ which gives
$$\nu \ge{ 1 \over (\l_{ikj}-1)\mu_{ij}(a_i,a_j)}+{\l_{ikj} \l_{jik} \over
\l_{ikj}-1} {(a_j,a_k) \over (a_j,a_i)(a_i, a_k)}.$$ If $k=j$ (and
so $l\ne j$) one similarly needs
$$\nu \ge {\l_{ijl} \over (\l_{ijl}-1)\mu_{ij}(a_i,a_j)}
+{\l_{jil} \over \l_{ijl}-1}{(a_j,a_l) \over
(a_j,a_i)(a_i,a_l)}.$$ \endpf

This lemma conveys the idea that if the set $\{(a_i,a_j)\}_{i\ne
j}$ is not ``too spread around'' then the group
$\G=\<T_{a_1},\cdots, T_{a_n}\>$ is free on $n$ generators, as
follows:

\begin{theorem}\label{3.2} Let $a_1,\cdots,a_n$ be $n\ge 3$ simple closed curves on a
surface $S$ such that $M \le m^2/6$ where
$M=\max\{(a_i,a_j)\}_{i\ne j}$ and $m=\min \{(a_i,a_j)\}_{i\ne
j}$. Then  $$\G=\<T_{a_1},\cdots, T_{a_n}\> \cong \mathbb F_n.$$
 More
generally, suppose that for all $i \ne j \ne k$ we have
$${(a_i,a_k) \over (a_i, a_j)(a_j,a_k)} \le {1 \over 6}.$$
Then the same conclusion holds.
\end{theorem}

\proof Put $\mu_{ij}=1$ and $\l_{ijk}=2$ in Lemma~\ref{3.1}. By
assumption, for all $i\ne j \ne k$,
$${(a_i,a_k) \over (a_i, a_j)(a_j,a_k)} \le {1 \over 6}.$$
This implies $(a_i,a_j) \ge 6$ for all $i\ne j$, since otherwise
it is impossible for both of
$${(a_i,a_k) \over (a_i, a_j)(a_j,a_k)} \  \textrm{and} \
{(a_j,a_k) \over (a_i, a_j)(a_i,a_k)}$$ to be $\le 1/6$.
Therefore, it is easily seen that $\nu=1$ satisfies the
requirements of Lemma~\ref{3.1}. \endpf

\setlength{\unitlength}{0.045in}

 \begin{figure}[ht!]\small
 $$\begin{picture}(0,0)(0,11)
   %\multiput(0,0)(5,0){25}{\line(0,1){150}}
   %\multiput(0,0)(0,5){25}{\line(1,0){150}}
   \put(-12,-7){$b$}
   \put(-4,-13){$a_2$}
  \put(1,10){$a_1$}
  \put(-23,-5){$\delta_1$}
  \put(22,-5){$\delta_2$}
   \end{picture}$$

  {\centerline{ \includegraphics[width=2.25in]{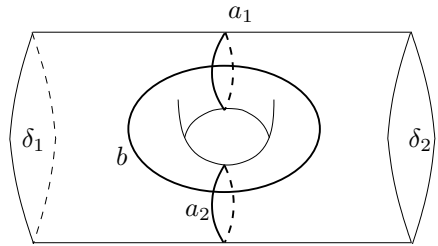}}}
     \caption{$(T_{a_1}T_{a_2}T_b)^4=T_{\delta_1}T_{\delta_2}$}
     \label{s22}
     \end{figure}

\section{Questions}

We end this paper by looking at some questions. Consider the group
$\G=\< T_{a_1}^{m_1}T_{a_2}^{m_2}, T_{b}^n \>$, where the simple
closed curves are defined in Figure~\ref{s22}, and they satisfy
the torus relation
$(T_{a_1}T_{a_2}T_b)^4=T_{\delta_1}T_{\delta_2}$. It is
interesting to find out if there are any torus-like relations.
Theorems~\ref{2multifree} and ~\ref{2multipa} will restrict the
search. In particular:

\rk{Question 1}Is it true that $\<T_{a_1}^2T_{a_2}, T_b
\>=\mathbb F_2$?

\rk{Question 2}Under what conditions is $\G=
\< T_{a_1}, \cdots, T_{a_n}\>$ relatively pseudo-Anosov?

\rk{Acknowledgments}The author was
partially supported by PSC-CUNY Research Grant 63463 00 32. He
thanks Marty Scharlemann, Darren Long and Daryl Copper for helpful
conversations and support while part of this work was being
completed. Thanks go to John McCarthy for posing the problem.
Thanks also go to the referee, for pointing out a few mistakes in
the proofs (which were subsequently fixed). The author would like
to thank Nikolai Ivanov and Benson Farb  for organizing a
wonderfully stimulating session on mapping class groups. Many
thanks to Dan Margalit for thoroughly reading the manuscript and
making a lot of constructive suggestions.
 Finally
thanks to Tara Brendle Owens for many encouraging and supporting
remarks.

 \Addresses\recd

 \end{document}